\documentclass[final]{siamltex}

\usepackage[normalem]{ulem}
\usepackage{amsmath}
\usepackage{amsfonts}
\usepackage{amssymb}
\usepackage{graphicx}
\usepackage{epstopdf,euscript}
\usepackage{url}
\usepackage{multirow}
\usepackage{array}
\usepackage{color}
\usepackage[ruled, linesnumbered]{algorithm2e}
\usepackage{mathtools}
\usepackage{subcaption}

\def\bq{\begin{quotation}}
\def\eq{\end{quotation}}

\def\n{\nu}

\def\2nm#1{\|#1\|_2}

\def\Ra#1{\mathrm{Range}(#1)}

\newcommand{\ars}[1]{\left[ \begin{array}{#1}}
\newcommand{\are}{\end{array} \right] }
\newcommand{\oars}[1]{\begin{array}{#1}}
\newcommand{\oare}{\end{array}}
\newcommand{\rars}[1]{\left( \begin{array}{#1}}
\newcommand{\rare}{\end{array} \right) }

\newcommand{\eqs}{\begin{eqnarray}}
\newcommand{\eqe}{\end{eqnarray}}
\newcommand{\eqsn}{\begin{eqnarray*}}
\newcommand{\eqen}{\end{eqnarray*}}

\newcommand{\bmp}[2]{\begin{minipage}#1{#2}}
\newcommand{\emp}{\end{minipage}}

\newcommand{\ens}{\begin{enumerate}}
\newcommand{\ene}{\end{enumerate}}

\newcommand{\its}{\begin{itemize}}
\newcommand{\ite}{\end{itemize}}

\newcommand{\des}{\begin{description}}
\newcommand{\dee}{\end{description}}

\def\defs{\begin{definition}}
\def\defe{\end{definition}}
\def\teos{\begin{theorem}}
\def\teoe{\end{theorem}}
\def\prfs{\begin{proof}}
\def\prfe{\end{proof}}
\def\exas{\begin{exampl}}
\def\exae{\end{exampl}}
\def\excs{\begin{exercise}}
\def\exce{\end{exercise}}
\def\cors{\begin{corollary}}
\def\core{\end{corollary}}

\newcommand{\ran}[1]{\mathrm{range}(#1)}

\newcommand{\cbfm}{\mbox{\boldmath${\mathit{m}}$} }

\newcommand{\bfsfp}{\mbox{\boldmath$\mathsf{p}$} }

\newcommand{\comment}[1]{} % makes its argument disappear

\newcommand{\bfPsi}{\boldsymbol{\Psi}}

\newcommand{\bfx}{{\bf x}}
\newcommand{\bfK}{{\bf K}}

\newcommand{\bfy}{{\bf y}}

\newcommand{\bfu}{{\bf u}}
\newcommand{\bfv}{{\bf v}}
\newcommand{\bfb}{{\bf b}}

\newcommand{\bfr}{{\bf r}}

\newcommand{\bfz}{{\bf z}}
\newcommand{\bfq}{{\bf q}}

\newcommand{\bfS}{{\bf S}}

\newcommand{\bfQ}{{\bf Q}}
\newcommand{\bfA}{{\bf A}}
\newcommand{\bfB}{{\bf B}}
\newcommand{\bfC}{{\bf C}}
\newcommand{\bfD}{{\bf D}}
\newcommand{\bfE}{{\bf E}}
\newcommand{\bfG}{{\bf G}}
\newcommand{\bfH}{{\bf H}}
\newcommand{\bfI}{{\bf I}}
\newcommand{\bfF}{{\bf F}}
\newcommand{\bfR}{{\bf R}}

\newcommand{\Co}{{\bf K}}

\newcommand{\bfX}{{\bf X}}

\newcommand{\bfY}{\widehat{\bf y}}
\newcommand{\bfYY}{{\bf Y}}
\newcommand{\bfU}{\widehat{\bf u}}
\newcommand{\bfUU}{{\bf U}}
\newcommand{\bfM}{\widehat{\bf m}}
\newcommand{\bfL}{{\bf \Lambda}}
\newcommand{\bfLL}{{\bf L}}
\newcommand{\bfN}{{\bf N}}
\newcommand{\bfZ}{{\bf Z}}
\newcommand{\bfV}{{\bf V}}
\newcommand{\bfW}{{\bf W}}

\newcommand{\bfg}{{\bf g}}
\newcommand{\bfe}{{\bf e}}

\newcommand{\IR}{{\mathbb{R}}}
\newcommand{\IC}{{\mathbb{C}}}

\newcommand{\bea}{\left[ \begin{array} }
\newcommand{\eea}{ \end{array} \right] }

\newcommand{\mrg}{\mbox{Range}}

\newcommand{\bDelta}{{{\Delta}}}

\title{Computing Reduced Order Models via Inner-Outer Krylov Recycling in Diffuse Optical Tomography\footnotemark[1]}

\author{ 
Meghan O'Connell\footnotemark[3]\ ,
Misha E. Kilmer\footnotemark[3]\ ,
Eric de Sturler\footnotemark[2]\ , and
Serkan Gugercin\footnotemark[2]
}

\begin{document}
\maketitle
\renewcommand{\thefootnote}{\fnsymbol{footnote}}
\footnotetext[1]{This material is based upon work supported by the National Science
Foundation under Grants No. {NSF-DMS} 1025327, {NSF DMS} 1217156 and 1217161,
{NSF-DMS} 0645347, and {NIH} R01-CA154774.}
\footnotetext[2]{Department of Mathematics, Virginia Tech, Blacksburg, VA 24061.}
\footnotetext[3]{Department of Mathematics, Tufts University, Medford, MA 02115.}
\renewcommand{\thefootnote}{\arabic{footnote}}

%%%%%%%%%%%%%%%ABSTRACT %%%%%%%%%%%%%%%%%
\begin{abstract}
  In nonlinear imaging problems whose forward model is described by a partial differential equation (PDE), the main computational bottleneck in solving the inverse problem is the need to solve many large-scale discretized PDEs 
  at each step of the optimization process.  In the context of absorption imaging in diffuse optical tomography, one approach to addressing this bottleneck proposed recently (de Sturler, et al, 2015) reformulates the viewing of the forward problem as a differential algebraic system, and then  employs  model order reduction (MOR).   However, the construction of the reduced model requires the solution of several full order problems (i.e. the full discretized PDE for multiple right-hand sides) to generate a candidate global basis.  This step is then followed by a rank-revealing factorization of the matrix containing the candidate basis in order to compress the basis to a size suitable
  for constructing the reduced transfer function.  The present paper addresses the costs associated with the global
  basis approximation in two ways.  First, we use the structure of the matrix to rewrite the full order transfer function, and corresponding derivatives, such that the full order systems to be solved are symmetric (positive definite in the zero frequency case).  Then we apply MOR to the new formulation of the problem. Second, we give an approach to 
  computing the global basis approximation dynamically as the full order systems are solved.  In this phase, only the incrementally new, relevant information is added to the existing global basis, and redundant information is not 
  computed.  This new approach is achieved by an inner-outer Krylov recycling approach which has potential use in other applications as well.   We show the 
  value of the new approach to approximate global basis computation on two DOT absorption image reconstruction problems.   
\end{abstract}

%%%%%%%%%%%%%%%%%%%% KEYWORDS %x%%%%%%%%%%%
%%%%%%%%%%%%%%%% AMS NUMBERS %%%%%%%%%%%%

%%%%%%%%%%%%%%% PAGE STYLE AND LABELING %%%%%%%%%%%%%
%\thispagestyle{plain}
\markboth{Kilmer, O'Connell, de Sturler, Gugercin}{Computing ROMs in DOT}

\section{Introduction}
Nonlinear inverse problems are becoming ubiquitous but remain very expensive to
solve, including medical image reconstruction, identification of
anomalous regions such as unexploded ordinance, land mines, and contaminant
plumes in the subsurface. In inverse problems we recover an image of an unknown
quantity of interest inside a given medium, such as the light absorption coefficient in
human tissue, using a mathematical model, the {\it forward model}, that relates
the image of the unknown quantity to the measured data. In this paper, we consider the
problem of imaging in diffuse optical tomography (DOT), where the
forward model is a  large-scale, discretized, partial differential equation describing
the photon fluence/flux through tissue. 
Other imaging problems (e.g. electrical impedance/resistance tomography, hydraulic tomography) are similar in that the forward problems are described by 
PDEs.     
  
The need to solve these large-scale, forward problems many times to recover
two or three-dimensional images represents the largest
computational impediment to effective, practical use of DOT.  
To this end, the use of reduced order modeling in the context of DOT imaging was considered in \cite{Gugetal2015}.  
The key observation in that work is that 
function evaluations for the underlying optimization problem (which requires solving for the parameters defining the parametric level set 
image)
may be viewed as transfer function evaluations along the imaginary axis, which motivates the use of
system-theoretic model order reduction methods. Specifically,
interpolatory parametric reduced models as surrogates for the full
forward model were used. In the DOT setting, such surrogate models were found to approximate
both the cost functional and the associated Jacobian with very little loss of
accuracy while at the same time drastically reducing the cost of the overall inversion process. 

Nevertheless, difficulties remain.  In order to determine the global basis to be used for the interpolatory projection,
several large-scale linear systems must still be solved.   Because of the structure and size of the systems,
the solves are typically done iteratively \cite{beattie2010isi}.  In recent work \cite{AhujaSturlerFengBenner2015,FenBen09a,FengBennerKorvink2013}, the authors investigate the use of Krylov subspace recycling for these systems to generate the model-order reduction (MOR) basis.   Recycling for shape based inversion for DOT was
investigated in \cite{KdS2006}, but the goal was solving the sequence of systems in the inversion
process, rather than for use in computing global basis matrices for use with MOR.     

In the present work, we present several computational advances for applying MOR to this inverse problem, but note that
{\it the results are applicable in a more general interpolatory MOR setting}.   The paper is organized as follows.
In Section \ref{sec:background}, we give the system theoretic notation, introduction to DOT, and relevant
information on computing reduced order models.   
In Section \ref{sec:symmetry}, we show that in a particular geometry and discretization, the transfer function
may be reformulated as a transfer function for a slightly smaller symmetric problem.   In Section \ref{sec:recycle}, we consider innovations in
leveraging Krylov recycling techniques to reduce the total amount of computation for the reduced global basis
matrix by eliminating redundant information dynamically.   An analysis is given in Section \ref{sec:analysis}.  Numerical
results are given in Section \ref{sec:numerical} and conclusions and future work are outlined in Section \ref{sec:conclusions}.

\section{Background}  \label{sec:background}

We begin by describing the forward model for DOT.  Then, we express the forward problem in the systems theoretic notation that we will use throughout the paper.  The first two sections 
borrow heavily from the presentation in our previous paper \cite{Gugetal2015}.  The third subsection gives background on the generation of the reduced transfer function, in preparation for the presentation of our new techniques for computing the reduced model global projection bases given in Sections
\ref{sec:symmetry} and \ref{sec:recycle}.  

\subsection{The DOT Problem}
Our image domain is a rectangular slab,
$\Omega = [a_1, b_1]\times[a_2, b_2]\times[a_3,b_3]$,
with the top ($x_3 = b_3$) and bottom ($x_3 = a_3$) surfaces denoted
by $\partial\Omega_{+}$ and $\partial\Omega_{-}$, respectively.   We use a diffusion model for the photon
flux/fluence \cite{arridge} $\eta(\bfx,t)$ driven by an input light source $g(\bfx,t)$.   The input light source is
one out of a set of $n_{src}$ possible sources
that are each physically stationary.  This means there are functions,
$b_j(\bfx)$, $j=1,\,\ldots,\, n_{src}$ such that $g(\bfx,t)=b_j(\bfx)u_j(t)$ for selected $j$.
Additionally the observations are made with a limited number of $n_{det}$ detectors observations.  We use $m_i(t)$ to denote the presumed stationary observations located on the bottom surface.
Using these definitions, the model for the diffusion and absorption of light is described via
\begin{align}
  \frac{1}{\nu}\frac{\partial}{\partial t}\eta(\bfx,t) & =
    \nabla \cdot \left(\, D(\bfx) \nabla \eta(\bfx,t) \,\right)
    - \mu(\bfx) \eta(\bfx,t) + b_j(\bfx)u_j(t), \quad \mbox{for }  \bfx\in \Omega ,
\label{eq:ModelPDE1} \\
    0 & = \eta(\bfx,t) +  2\,{\cal A}\,D(\bfx)\,
      \frac{\partial}{\partial \xi}\eta(\bfx,t) , \quad \mbox{for }
      \bfx\in \partial\Omega_{\pm} ,
\label{eq:ModelPDE2} \\
    0 & = \eta(\bfx,t), \quad \mbox{for } x_1 = a_1 \;\mathrm{or}\; x_1 = b_1
      \;\mathrm{or}\; x_2 = a_2 \;\mathrm{or}\; x_2 = b_2 ,
\label{eq:ModelPDE3} \\
    m_i(t) & =  \int_{\partial\Omega} c_i(\bfx)\eta(\bfx,t)\,d\bfx\quad
      \mbox{ for }i=1,\,\ldots,\, n_{det}
\label{eq:ModelPDE4}
\end{align}
(see \cite[p. R56]{arridge}).  Here, the vector $\bfx=(x_1,x_2,x_3)^T$
refers to spatial location, ${\cal A}$ is a constant defining the particular diffusive boundary reflection (see \cite[p.R50]{arridge}), and $D(\bfx)$ and $\mu(\bfx)$ denote diffusion and absorption
coefficients, respectively.  Also,
$\xi$ denotes the outward unit normal and $\nu$ is the speed of
light in the medium.

%The scalar fields defined by $D(\bfx)$ and $\mu(\bfx)$ are assumed unknown or partially
%known.  
The inverse problem consists of utilizing observations, $\cbfm(t)$, made when the system is
illuminated by a variety of source signals, $\bfu(t)$,
to more accurately determine $D(\bfx)$ and $\mu(\bfx)$.
For purposes of this work, we assume the diffusivity $D(\bfx)$ is
known (a common assumption in DOT breast tissue imaging) 
and that only the absorption field, $\mu(\bfx)$, must be recovered.
We also assume that the absorption field, $\mu(\cdot)$, although unknown, is expressible
in terms of a finite set of parameters, $\bfsfp=[p_1,\ldots,\,p_{\ell}]^T$.
%An effective parametrization of $\mu(\cdot)$ is crucial.   
We will assume that parametric level sets (PaLS), developed in \cite{Aghasi_etal11} and used in the context
of DOT imaging in \cite{Gugetal2015}, have been used for $\mu(\cdot)=\mu(\cdot,\bfsfp)$.
The inverse problem is often solved by moving to the frequency domain and using frequency domain data (e.g. data using frequency modulated light).  However, we adopt the approach in \cite{Gugetal2015} of reformulating the forward problem in dynamical systems notation before moving to the frequency domain, so that 
we can see how to employ model reduction in the context of solving the inverse problem.

%One may use FEM or FD to
The discretization of  (\ref{eq:ModelPDE1})-(\ref{eq:ModelPDE4}) can be done by finite element or finite difference techniques.
%finite element methods and finite difference
%methods among them.
This gives the following differential-algebraic system,
\begin{equation}  \label{processModelDynSys}
  \frac{1}{\nu}\bfE\,  \dot{\bfy}(t;\bfsfp)  =-\bfA(\bfsfp)\bfy(t;\bfsfp) +
    \bfB\bfu(t)\quad\mbox{with}\quad
    \cbfm(t;\bfsfp) =\bfC^T\bfy(t;\bfsfp)
\end{equation}
where $\bfy$ denotes the discretized photon flux,
$\cbfm=[m_1,\,\ldots,\, m_{n_{det}}]^T$ is the vector of detector outputs, $\bfC^T\bfy$ constitutes a
set of quadrature rules with the approximate
photon flux; the columns of $\bfB$ are discretizations of the
source ``footprints" $b_j(\bfx)$ for $j=1,\,\ldots,\, n_{src}$;
$\bfA(\bfsfp)=\bfA_{0}+\bfA_{1}(\bfsfp)$ with
$\bfA_{0}$ and $\bfA_{1}(\bfsfp)$ discretizations of the diffusion and absorption
terms, respectively ($\bfA_{1}(\bfsfp)$ inherits the
absorption field parametrization, $\mu(\cdot,\bfsfp)$).
$\bfE$ is singular due to the inclusion of the discretized Robin condition (\ref{eq:ModelPDE2}) as an algebraic constraint.  This fact will become important in Section \ref{sec:symmetry}. 

Let $\bfY(\omega; \bfsfp)$, $\bfU(\omega)$, and $\bfM(\omega; \bfsfp)$ denote
the Fourier transforms of $\bfy(t;\bfsfp)$, $\bfu(t)$, and $\cbfm(t;\bfsfp)$, respectively.
Taking the Fourier transform of (\ref{processModelDynSys}) and rearranging, we get
\begin{equation} \label{FullOrdTransFnc}
  \bfM(\omega;\bfsfp) =\bfPsi\!\left(\omega;\bfsfp\right)\,\bfU(\omega)\quad \mbox{where}\quad
  \bfPsi(\omega;\bfsfp)=\bfC^T \left(\frac{\imath\;\!\omega}{\nu}\,\bfE\, +\bfA(\bfsfp)\right)^{-1}\bfB,
\end{equation}
where $\omega \in \IR$, and $\bfPsi(\omega;\bfsfp)$  is known as the \emph{frequency response}
of the dynamical system defined in (\ref{processModelDynSys}), though we will often refer to it as the transfer function\footnote{In describing linear dynamical systems, usually the {\it transfer function}
$\bfPsi(s;\bfsfp)=\bfC^T \left(\frac{s}{\nu}\,\bfE\, +\bfA(\bfsfp)\right)^{-1}\bfB$ is used where  $s\in\IC$ and is not restricted to the imaginary axis. Here, though, the measurements are made only on the imaginary axis and it is enough to take $s = \imath\;\! \omega$ with $\omega \in \IR$.}.

For any absorption field, $\mu(\cdot,\bfsfp)$, associated with $\bfsfp$,  the vector of
(estimated) observations for the
$i^{\rm th}$ input source at frequency $\omega_j$, as predicted by the forward model in
the frequency domain, will be
denoted as $ \bfM_i(\omega_j;\bfsfp)\in\mathbb{C}^{n_{det}}$.
Stacking the predicted observation vectors for all $n_{src}$ sources and $n_{\omega}$ frequencies, we obtain
\eqs \label{eq:TotalObserv}
{\cal M}(\bfsfp) = [  \bfM_1(\omega_1;\bfsfp)^T,\, \ldots,\,\bfM_1(\omega_{n_\omega};\bfsfp)^T,\,\bfM_2(\omega_1;\bfsfp)^T, \ldots,
\bfM_{n_{src}}(\omega_{n_\omega};\bfsfp)^T ]^T,~~~
\eqe
which is a (complex) vector of dimension $n_{det}\cdot n_{src} \cdot n_{\omega} $.    We construct the corresponding empirical data vector, $\mathbb{D}$,
from acquired data.
The optimization problem that must be solved is
\begin{equation} \label{eq:InvProb} \min_{\bfsfp \in \mathbb{R}^\ell} \| {\cal M}(\bfsfp) - \mathbb{D} \|_2. \end{equation}
We note that mathematically, what was just described is equivalent to the usual approach for DOT imaging using frequency modulated data.  What is different is using the dynamical systems interpretation, as was developed in \cite{Gugetal2015}, so that we can leverage efficiencies from a systems theoretic perspective on solving the
parametric inversion problem.   

\subsection{A systems theoretic perspective on parametric inversion}

From (\ref{FullOrdTransFnc}) and (\ref{eq:TotalObserv}), it follows that a single evaluation of ${\cal M}(\bfsfp) - \mathbb{D}$
involves computing (for all $i$ and $j$)
\eqs \label{eq:CYi_om_j}
  \bfM_i(\omega_j; \bfsfp)=
    \bfPsi( \omega_j) \bfU_i(\omega_j) , \quad
\eqe
where $\bfPsi(\omega;\bfsfp)$ is the frequency response defined in (\ref{FullOrdTransFnc})
and $ \bfU_i(\omega_j) = e^{\imath \omega_j}\mathbf{e}_i$ where $\mathbf{e}_i$ is the $i^{\rm th}$
column of the identity matrix; i.e., $\bfU_i(\omega_j)$ excites the $i^{\rm th}$ source location.
The key observation is that as a result, an objective function evaluation at parameter vector $\bfsfp_k$ requires solving the block linear systems
\begin{equation} \label{eq:fwd} \left(\frac{\imath\;\! \omega_j}{\nu}\,\bfE\, +\bfA(\bfsfp_k)\right)
\mathbf{Y}_{k,j} = \bfB\,  \qquad \bfB \in \mathbb{R}^{n \times n_{src}}, \qquad j=1,\ldots,n_\omega .\end{equation}
   Here, 
\[ \mathbf{Y}_{k,j} := [\bfY_1(\omega_j,\bfsfp_k),\ldots,\bfY_{n_s}(\omega_j,\bfsfp_k) ] .\]
It is important to note that each column of the matrix $\bfB$ is, in our application, a multiple of the $i^{\rm th}$ column of the identity matrix, corresponding to the $i^{\rm th}$ source location.

To solve the nonlinear inverse problem (\ref{eq:InvProb}) for the parameters, the
 Jacobian is constructed using an adjoint-type (or co-state)
approach that exploits the fact that the number of detectors
is roughly equal to the number of sources,
as discussed in \cite{HabAsch00} and \cite[p. 88]{Vogel}.
Using (\ref{FullOrdTransFnc}) and (\ref{eq:CYi_om_j})
and differentiating $ \bfM_i(\omega_j;\bfsfp)$
with respect to the $k$th component of $\bfsfp$,
\begin{equation} \label{eq:jaccalc}
  \frac{\partial }{\partial p_k} \bfM_i(\omega_j;\bfsfp)  =  \frac{\partial }{\partial
    p_k}\left[\bfPsi\!\left(\omega_j;\bfsfp\right)\right]\,\bfU_i(\omega_j)
    =  - \bfZ(\omega_j;\bfsfp)^T\,
    \frac{\partial }{\partial p_k}\bfA(\bfsfp)\, \bfY_i(\omega_j; \bfsfp),
\end{equation}
where, for each $\omega_j$ and any $\bfsfp$ we can compute $\bfZ(\omega_j;\bfsfp)$ from
\begin{equation} \label{eq:adj}
\left(\frac{\imath\;\! \omega_j}{\nu}\,\bfE\, +\bfA(\bfsfp)\right)^T \bfZ(\omega_j;\bfsfp) =  \bfC, \qquad 
\bfC \in \mathbb{R}^{n \times n_{src}}.
\end{equation}
Here, the columns of $\bfC$ also correspond to columns of the identity matrix, with a non-zero entry appearing at the $j$th detector location.

The matrices $\frac{\partial }{\partial p_k}\bfA(\bfsfp)$ need to be computed
only once for all parameters. 
Except for the cost of computing $\bfY_i(\omega_j; \bfsfp)$ (which was completed already for the function evaluation), the computational cost of evaluating the Jacobian at a parameter vector, for all $\omega_j$,
consists mainly of the cost of computing the solution to the block systems (\ref{eq:adj}).

In sum, the critical bottleneck in solving the inverse problem
is associated with repeatedly solving (\ref{eq:fwd}) and (\ref{eq:adj}).

\subsection{Global Basis Computation}  \label{ssec:globalcomp}
As done in \cite{Gugetal2015} for the DOT problem, we therefore seek a surrogate 
function ${\bfPsi}_r(\omega;\bfsfp)$
that is much cheaper to evaluate, yet provides a high-fidelity approximation to
$\bfPsi(\omega;\bfsfp)$  {\it over parameters and  frequencies of interest}.  The frequencies
are dictated to us by the experimental set-up.   The parameters of interest are those which would be chosen
by the optimization routine if it were run using the full order model.  
Likewise, we require that $\nabla_{\bfsfp}{\bfPsi}_r(\omega;\bfsfp)$ is much cheaper
to evaluate and that
$\nabla_{\bfsfp}\bfPsi(\omega;\bfsfp)\approx \nabla_{\bfsfp}{\bfPsi}_r(\omega;\bfsfp)$
over the same range of arguments.

The surrogate parametric model is obtained using projection (see \cite{Benner2015survey,Gugetal2015}).   
Suppose $\bfA \in \mathbb{R}^{n \times n}$ and full rank matrices $\bfV \in \IC^{n \times r}$ and $\bfW\in \IC^{n \times r}$
are given. Assuming the full state $\bfy(t;\bfsfp)$ evolves near the $r$-dimensional
subspace $\Ra{\bfV}$, one has
$\bfy(t;\bfsfp) \approx \bfV \widehat{\bfy}(t;\bfsfp)$.  If we enforce a Petrov-Galerkin condition
we obtain a reduced system with
the reduced matrices given by
\begin{align} \label{eq:PGProj}
  {\bfE_r} = \bfW^{T}\bfE\bfV, \quad
    {\bfA_r}(\bfsfp) = \bfW^{T} \bfA(\bfsfp) \bfV,
    \quad {\bfB_r}= \bfW^{T}\bfB, \quad \mbox{ and } \quad {\bfC_r} = \bfV^T \bfC.
\end{align}
The reduced transfer function is then ${\bfPsi}_r = {\bfC}_r^T ( \frac{\imath\;\! \omega}{\nu} \bfE_r + {\bfA}_r(\bfsfp))^{-1}{\bfB_r}$.   The goal is to choose $\bfV, \bfW$ so that the transfer function evaluated at $\bfsfp, \omega$
(likewise the Jacobian) will match (or be close to) the reduced transfer function evaluated at those parameters
and frequencies of interest.      

%In \cite{dSetal2015}, a ``global basis approach'' is used to construct
%matrices $\bfV$ and  $\bfW$  that attempt to cover the parametric information for the
%entire parameter space by including information from multiple parameter values.
%Suppose $K$ parametric sample points  
%$\boldsymbol{\pi}_1,\ldots,\boldsymbol{\pi}_K$, for $i=1,\ldots,K$,
%have been determined.  
Referring to (\ref{eq:fwd}), we let
\begin{equation}   \label{eq:localVi}  % \label{eq:Vij}
  %\bfV_i  =  \left[ \bfV_{i,1}, \bfV_{i,2}, \ldots, \bfV_{i,n_\o}\right], \mbox{ and then }  
  {\bfYY} := [\bfYY_{1,1},\ldots,\bfYY_{1,n_\omega},\bfYY_{2,1},\ldots,\bfYY_{2,n_\omega},\ldots,
   \bfYY_{K,n_\omega}] . 
  \end{equation}
Then, ideally, we define $\bfV$ to correspond to the left singular vectors corresponding to the non-zero singular values of the matrix ${\bfYY}$.
Similar steps are applied to construct $\bfW$ from $\bfZ_{i,j}$ via (\ref{eq:adj}) for $i=1,\ldots,K$
and $j=1,\ldots,n_\omega$.   Taking the exact left singular vectors for non-zero singular values
to define both $\bfV, \bfW$ ensures \cite{BauBBG09} that the corresponding reduced transfer function
will match the transfer function evaluation at {\em every}
$(\omega,\bfsfp) =  (\omega_j,\boldsymbol{\pi}_i)$ for
$j=1,\ldots,n_\omega$ and $i=1,\ldots,K$.   A similar result holds for the respective derivative computations.

%However, in practice, the singular values of ${\bfYY}$ fall off several orders of magnitude.   Thus, we truncate the left singular vectors of ${\bfYY}$ according
%to a threshold and call this $\bfV$, and likewise for $\bfW$.  Clearly, we have a more compressed basis and less work to form the reduced order model, but we give up the 
%exact matching for approximate matching.   

The one-sided global basis approach, which is used in the context of DOT in \cite{Gugetal2015}, refers to using $\bfV \leftarrow [\bfV, \bfW]$ and $\bfW \leftarrow [\bfV,\bfW]$ in (\ref{eq:PGProj}).  
%In other words, both the forward and adjoint information is computed, and SVDs are used on the consolidated
%information, then the kept left singular vectors are concatenated into a single $\bfV$.  
Clearly,
if $\bfA(\bfsfp), \bfE$ are symmetric (Hermitian), the reduced counterparts defined in 
(\ref{eq:PGProj}) are also symmetric (Hermitian).   For the use of model reduction in other optimization  and inverse problem applications, we refer
the reader to
\cite{Arian2002,hinze2005proper,Kunisch2008,Antil2011,Antil2012,Druskin2011solution,bashir08,borcea2012model,
yue2013} and the references therein.

In the preceding description, we have assumed (as elsewhere in the literature), that 
several full order problems for the $\bfYY_{ij}$ ($ \bfZ_{ij}$) have been solved prior to selecting the global basis.  
 In the context of solving our inverse problem, however, we do not a priori pick parameter values, but rather use the parameter sequence defined at the start of the optimization to construct $\bfV$ from the required full order model (FOM) solves. %in the first few steps of the optimization process, and
%we will describe this new approach in Section \ref{sec:recycle}. 
It was observed in \cite{Gugetal2015} for our DOT application that then the singular values of $\bfYY$ ($\bfZ$) fall off by several orders of magnitude.  Thus, truncated SVDs are used to obtain the $\bfV, \bfW$ which are subsequently concatenated to form one large projection matrix $\bfV$.  
This means that in first solving for all the $\bfYY_{ij}, \bfZ_{ij}$, we must have computed redundant 
information that we subsequently must squeeze out through via SVDs on the concatenated matrices.   The purpose of this paper is to construct an approximate symmetric global basis matrix {\em without first solving each linear system in \ref{eq:fwd} and \ref{eq:adj}} so that we can also avoid the additional post processing step of computing truncated SVDs. 
  
First, however, we show in \S \ref{sec:symmetry} how we can capitalize on the structure of the system matrix in our case to 
design a ROM scheme on a slightly different system matrix/variables.   Then in Section \ref{sec:recycle} we exploit this property to develop a
particularly efficient algorithm and corresponding analysis for generating columns of the projection matrix $\bfV$.

%%%%%%%%%%%%%%%%%%%%%%%%%%%%%%%%%%%%%%%%%%%%
\section{Rewriting Transfer Function and Derivatives}  \label{sec:symmetry}

In this section, we show that the structure of the linear system allows us to express the 
transfer function in terms of a symmetric positive definite matrix.
%which is differentthan the one in (\ref{FullOrdTransFnc}).   
%The value of this approach is that
%the new matrix will, under our discretization, be SPD.  
This leads to a more efficient method for generating the global basis, and is of benefit in some theoretical arguments.

We follow the same discretization as in \cite{KdS2006}, in which we use second order centered differences away from the boundary, and first order discretization to implement the Robin boundary condition.  The matrix $\frac{\imath\;\! \omega}{\nu}\bfE + \bfA(\bfsfp)$ appearing in the definition of the full
order transfer function (\ref{FullOrdTransFnc}) in the 2D case, has the following
block structure:
\begin{equation}  \label{eq:block1}
\bf \left[ \begin{array}{cc}
      \bfG   & \bfD_1 \\
      \bfD_2 & (\bfF(\bfsfp) + \frac{\imath\;\! \omega h^2}{\nu} \bfI)
    \end{array} \right],
\end{equation}
where we have ordered the $N_x N_y$ unknowns such that the boundary unknowns (lexicographically ordered) appear
first, followed by lexicographical ordering of internal points.

Furthermore, regarding the blocks in (\ref{eq:block1}),
\begin{itemize}
  \item $\bfG$ is an invertible diagonal matrix,
  \item $\bfD_1$ has at most one nonzero per row, and these occur only in the first
$N_x N_y$ and last $N_x N_y$ columns,
  \item $\bfD_2$, although it has different entries, has the same sparsity pattern as $\bfD_1^T$.
\end{itemize}

We note that the matrix $\bfA(\bfsfp)$ is not symmetric; however, as is shown in 
\cite{KdS2006},  the Schur complement (for the $\omega = 0$ case), which
is given by $\bfF(\bfsfp) - \bfD_2 \bfG^{-1} \bfD_1$, is symmetric and positive definite.
We will now
investigate why this fact is particularly useful with regard to specification of the transfer function.

   %%%%%%%%%%%%%%%%%%%%%%%%%%%%%%%%%%%%%%%%%%
\subsection{Transfer Function Revisited}

The columns of $\bfC$ and $\bfB$ are scaled columns from an $N_x N_y \times N_x N_y$ identity matrix.  Since sources
and detectors appear on the boundary, this means that partitioning conformably with $\frac{\imath\;\! \omega}{\nu} \bfE + \bfA(\bfsfp)$ as in (\ref{eq:block1}) we obtain
\[ \bfC = \bea{c} \bfC_1 \\  {\bf 0} \eea, \bfB = \bea{c} \bfB_1 \\ {\bf 0} \eea . \]
Since the sources and
detectors are not co-located, it follows that $\bfC_1^T \bfB_1 = 0$.

Let us assume for ease of exposition that $\omega = 0$ in (\ref{eq:block1}) above.
In this case, it is the inverse of $\bfA(\bfsfp)$ that appears in the definition of the transfer function. 
Let the matrix $\bfA(\bfsfp)^{-1}$ have the following block structure:
\[
\bf \left[ \begin{array}{cc}
      \bfH   & \bfS_1 \\
      \bfS_2 & \bfN
    \end{array} \right]. \]
From $\bfA(\bfsfp) \bfA(\bfsfp)^{-1} = \bfI$, we have the following three expressions that will be
helpful in rewriting the transfer function and in computing corresponding measurement derivatives:
\begin{equation} \label{eq:H} \bfH = [\bfG - \bfD_1^T \bfF^{-1} \bfD_2 ]^{-1},   \end{equation}   
\begin{equation} \label{eq:Sone} \bfS_1 = -\bfG^{-1}\bfD_1[\bfF - \bfD_2\bfG^{-1}\bfD_1]^{-1}, \end{equation}
\begin{equation} \label{eq:Stwo} \bfS_2 = -[\bfF - \bfD_2\bfG^{-1}\bfD_1]^{-1}\bfD_2\bfG^{-1}.\end{equation} 

It is straightforward to show that
\[ \bfPsi(0,\bfsfp) = \bfC^{T} ( \bfA(\bfsfp) )^{-1} \bfB = 
\bfC_1^T \bfH \bfB_1. \]
Now $\bfF(\bfsfp) = \bfLL + \mbox{diag}(\mu(\bfsfp))$, where $\bfLL$ is the discretization of the Laplacian at the internal nodes multiplied by the (constant) diffusion coefficient, and $\mbox{diag}(\mu(\bfsfp))$ is a non-negative diagonal matrix.  Thus $\bfF(\bfsfp)$ is SPD, so we express it in terms of 
its eigendecomposition,
$ \bfF = \bfQ \bfL \bfQ^T $, so that
\begin{eqnarray*}
 \bfH &=& [\bfG - \bfD_1^T \bfQ \bfL^{-1/2} \bfL^{-1/2} \bfQ^T \bfD_2 ]^{-1} \\
       &= & \bfG^{-1} + \bfG^{-1} \bfD_1^T[ \bfF - \bfD_2 \bfG^{-1} \bfD_1]^{-1} \bfD_2 \bfG^{-1}, \end{eqnarray*}
by the Sherman-Morrison-Woodbury formula.   

Since the sources and detectors are not co-located, and because $\bfG$ is diagonal, it follows that 
\[ \bfPsi(0,\bfsfp) = \underbrace{\bfC_1^T \bfG^{-1}\bfD_1^T}_{\tilde{\bfC}^T} \underbrace{[ \bfF - \bfD_2 \bfG^{-1} \bfD_1]^{-1}}_{\tilde{\bfA}(\bfsfp)^{-1}} \underbrace{\bfD_2 \bfG^{-1} \bfB_1}_{\tilde{\bfB}}. \]
Importantly, the structure of the matrices involved means that $\tilde{\bfC}$ and $\tilde{\bfB}$ maintain 
the same structure as $\bfC_1$ and $\bfB_1$: namely, they contain multiples of columns of the $N_x N_y$ 
identity matrix, so those matrices can be considered as `effective' sources and receivers.  

Using a similar argument, it is also straightforward to show that if $\omega \not= 0$, 
\[ \bfPsi(\omega, \bfsfp) = \tilde{\bfC}^T \left( \frac{\imath\;\! \omega}{\nu} \bfI + \tilde{\bfA}(\bfsfp) \right)^{-1} \tilde{\bfB}^T.\]
This means that the transfer function (for the 0 frequency case) is in fact expressed using an SPD matrix to represent the system matrix (complex symmetric if
$\omega$ is non-zero).   

\subsection{Systems Theoretic Interpretation}
 
 Even though the new expression in the last section was derived from matrix analysis, there is also 
 a systems theoretic interpretation that will lead to the same expression. 
 
 Specifically, using the same reordering of unknowns as in the previous section, we observe that the singular matrix $\bfE$ has the structure
 \[ \bea{cc} 0 & 0 \\ 0 & \bfI \eea, \]
 so that the system is a differential-algebraic system and, due to the structure of $\bfE$, of index 1.  If we partition the state vector conformably with 
 the matrices $\bfA(\bfsfp), \bfE$, and $\bfB$, then $\bfE \dot{\bfy} = \bfA(\bfsfp) \bfy + \bfB \bfu$ becomes
 \begin{eqnarray}  0 & = & \bfG \bfy_1 + \bfD_1 \bfy_2 \label{eq:first} \\
                          \dot{\bfy}_2 & = & \bfD_2 \bfy_1 + \bfF \bfy_2. \label{eq:second} \end{eqnarray}
                          If we solve for $\bfy_1$ in (\ref{eq:first}) and plug the result into (\ref{eq:second}) and rearrange, we
                          obtain
                          \[ \dot{\bfy}_2 = \tilde{\bfA}(\bfsfp)\bfy_2 - \tilde{\bfB} \bfu. \]
                          Now the same calculations in the previous subsection lead us to 
                          \[ {\bf m}(t) = \bfC^T \bfy = \bfC_1^T \bfy_1 = -\tilde{\bfC}^T\bfy_2.\]   
                          Posing the two previous equations in the Fourier domain yields the transfer function
                          \[ \Psi(s,\bfsfp) = \tilde{\bfC}^T (s \bfI - \tilde{\bfA}(\bfsfp))^{-1} \tilde{\bfB} .\]

This manipulation here corresponds to decomposing the transfer function of a DAE as the sum of the strictly proper part and polynomial part. A satisfactory reduced model should match the polynomial component exactly. In the interpolatory MOR setting, this means that the strictly proper part of the reduced transfer function interpolates that of the full-order model. For details on interpolatory MOR of DAEs, we refer the reader to \cite{gugercin2013model}.
For the DOT problem we consider here, this derivation illustrates that the transfer function of the index-1 DAE does not contain any polynomial part and is strictly proper. 

%%%%%%%%%%%%%%%%%%%%%%%%%%%%%%%%%%%%%%5555
\subsection{Derivative Computation Revisited} \label{sec:derivative}

Again using $\omega = 0$ for simplicity,
define the $n_{src} \times n_{det}$ matrix 
$$ M(0,\bfsfp) = [\bfM_1(0,\bfsfp),\ldots,\bfM_{n_{det}}(0,\bfsfp)].$$  
We use the ``vec'' command to map a matrix in $\mathbb{R}^{m_1 \times m_2}$ to a vector in $\mathbb{R}^{m_1 m_2}$ by unstacking the columns of the argument from left to right.
Note that the vector  ${\rm{vec}}(\frac{\partial}{\partial p_k} M) \in \mathbb{R}^{n_{src} n_{det}}$ gives the
$k^{\rm th}$ column of the Jacobian matrix.   
Using ( \ref{eq:fwd} - \ref{eq:adj}),  we can write 
\[ \frac{\partial }{\partial p_k} M(\bfsfp)= - \bfC^T\bfA(\bfsfp)^{-1}
    \frac{\partial }{\partial p_k}\bfA(\bfsfp)\bfA(\bfsfp)^{-1}\bfB. \]
This calculation can be carried out with only reference to $\tilde{\bfA}(\bfsfp)$, as follows. 

Since the boundary terms have no absorption, then under the same finite difference discretization scheme and
ordering of unknowns as before, 
\[ \frac{\partial }{\partial p_k}\bfA(\bfsfp) = \bea{cc} \bf 0 & \bf 0 \\ \bf 0 & \Delta \eea, \]
for a diagonal matrix $\Delta := \frac{\partial}{\partial p_k} \tilde{\bfA}(\bfsfp).$  
Thus,
\begin{eqnarray*}
\frac{\partial }{\partial p_k} M(0,\bfsfp) &=& - \begin{bmatrix} \bfC_1^T & \bf 0 \end{bmatrix} \bea{cc} \bfH & \bfS_1 \\ \bfS_2 & \bfN \eea \bea{cc} \bf 0 & \bf 0 \\ \bf 0 & \Delta \eea \bea{cc} \bfH & \bfS_1 \\ \bfS_2 & \bfN \eea \bea{cc} \bfB_1 \\ \bf 0 \eea\\
    &=& - \bfC_1^T\bfS_1\Delta\bfS_2\bfB_1. \end{eqnarray*}
Now using (\ref{eq:Sone}) and (\ref{eq:Stwo}), we have 
\begin{equation} \label{eq:newjac}
\frac{\partial }{\partial p_k} M(0,\bfsfp)= - \underbrace{\bfC_1^T \bfG^{-1}\bfD_1}_{\tilde{\bfC}^T} \underbrace{[ \bfF - \bfD_2 \bfG^{-1} \bfD_1]^{-1}}_{\tilde{\bfA}(\bfsfp)^{-1}} {\Delta} \underbrace{[ \bfF - \bfD_2 \bfG^{-1} \bfD_1]^{-1}}_{\tilde{\bfA}(\bfsfp)^{-1}}\underbrace{\bfD_2 \bfG^{-1} \bfB_1}_{\tilde{\bfB}}.
\end{equation}
This means that the necessary derivatives, for the 0 frequency case, can be computed from the same SPD matrix as the transfer function. 

%%%%%%%%%%%%%%%%%%%%%%%%%%%%%%%%%%%
\subsection{A New View on Generating the ROM}
We have expressed the original transfer function and corresponding derivatives in 
terms of an SPD (for $\omega = 0$) matrix of size $(N_y-2)N_x \times (N_y -2)N_x$.  So, we look for a ROM corresponding to this slightly smaller system matrix and its system variables.     

Following the discussion in Section \ref{sec:background} of the one-sided global basis projection approach, 
we need $\bfV \in \mathbb{C}^{n \times r}$ and we define\footnote{The $\bfV$ as we generate it will typically not
have orthonormal columns, hence the need to specify $\tilde{\bfE}_r$.} 
\begin{equation} \label{eq:newreduced} \tilde{\bfE}_r = \bfV^T \bfV, \qquad \tilde{\bfA}_r(\bfsfp) = \bfV^T \tilde{\bfA}(\bfsfp) \bfV,\qquad \tilde{\bfB}_r = \bfV^T \tilde{\bfB}, \qquad \tilde{\bfC}_r = \bfV^T \tilde{\bfC};\end{equation}
so, the reduced 
transfer function is $\tilde{\bfPsi}_r = \tilde{\bfC}_r^T (\frac{\imath\;\! \omega}{\nu} \tilde{\bfE}_r + \tilde{\bfA}_r(\bfsfp))^{-1} \tilde{\bfB}_r$.   

Since $\tilde{\bfA}$ is SPD we do not need to solve the forward and adjoint problems separately to generate
$\bfV$.  Instead, we can solve
\begin{equation} \label{eq:solves}
\left(\tilde{\bfA}(\bfsfp_k) + \frac{\imath\;\! \omega_j}{\nu} \bfI\right) \bfX_{k,j} = [\tilde{\bfB},\tilde{\bfC}]
\end{equation}
for appropriate choices of parameters $\bfsfp_k$, $k=1,\ldots,K$ and frequencies $\omega_j, j=1,\ldots,\n_{\omega}$.   

For the remainder of the paper, we will assume that only data for $\omega=0$ has been provided.  This helps with
the remaining exposition of our new approach.  Indeed, only results for the zero frequency case are provided in \cite{Gugetal2015}.  Moreover, treatment 
of non-zero $\omega$ is not a trivial extension of the algorithm provided, and in the interest of both paper focus and space, we relegate extensions for non-zero $\omega$ to forthcoming work. 

The objective in reduced order modeling is to create a surrogate transfer function, 
$\tilde{\bfPsi}_r(\omega;\bfsfp)$, that provides a high-fidelity approximation to
$\tilde{\bfPsi}(\omega;\bfsfp)$ as well as ensuring $\nabla_{\bfsfp}\tilde{\bfPsi}(\omega;\bfsfp)\approx \nabla_{\bfsfp}{\tilde{\bfPsi}}_r(\omega;\bfsfp)$.  
The following theorem, which follows from \cite{BauBBG09}, shows how to construct $\bfV$ to guarantee this for the symmetric DOT-PaLs in the zero frequency case.   
\teos
Suppose $\tilde{\bfA}(\bfsfp)$ is continuously differentiable in a neighborhood of $\hat{\bfsfp} \in \mathbb{R}^l$. Let both $\tilde{\bfA}(\hat{\bfsfp})$ and $\tilde{\bfA}_r(\hat{\bfsfp})$ be invertible. 
If $\tilde{\bfA}(\hat{\bfsfp})^{-1} \bfB$ and $\left( \bfC \tilde{\bfA}(\hat{\bfsfp})^{-1} \right) ^T$ are in $\ran{\bfV}$, then the reduced parametric model satisfies 
$\tilde{\bfPsi}(0;\hat{\bfsfp}) = \tilde{\bfPsi}_r(0;\hat{\bfsfp})$ and $\nabla_{\bfsfp}\tilde{\bfPsi}(0;\hat{\bfsfp}) = \nabla_{\bfsfp}\tilde{\bfPsi}_r(0;\hat{\bfsfp})$.
\teoe   
 
In accordance with the discussion at the end of \S \ref{ssec:globalcomp}, the approach to obtaining $\bfV$ for the solution to the DOT problem would consist of the following steps.  Note that what follows is the approach of \cite{Gugetal2015}, but modified
for our newly formulated transfer function representation. 
   
%\noindent{Algorithm 1:  Generate Global Basis} 
\begin{algorithm}[ht]
   \begin{enumerate}
      \item Solve the systems (\ref{eq:solves}) for $k=1,\ldots,K_*$ for the first $\bfsfp_k$ parameter vectors
      produced by the optimization;
      \item Concatenate the block solutions into a large block matrix
      \[ [\bfX_{1,0},\bfX_{2,0},\ldots,\bfX_{K,0} ].\]
      \item Set $\bfV$ to be the matrix of the first $r$ left singular vectors of the above matrix.  This gives
      a reduced order model, (\ref{eq:newreduced}), of dimension $r$.  
      \end{enumerate}
      \caption{Generate Symmetric Global Basis via Truncated SVD}
      \label{alg:alg1}
\end{algorithm}
      
      \medskip
 On the one hand, we need to compute each $\bfX_{k,0}$, because we need this to compute function and 
 Jacobian evaluations at steps $1$ through $K_*$ of the optimization problem.   Clearly, as these are large and sparse systems with SPD matrices, we should employ a Krylov subspace algorithm to solve the individual systems.  
 On the other hand, for suitable $K_*$, some of the
 information that we generate and put into the concatenated matrix is redundant (or nearly so), as evidenced
 by the rapid decay of the singular values of the concatenated matrix.   This means in terms of generating a
 global basis matrix $\bfV$, we have computed information that we don't really need (adding to the cost) and 
 we thus incur the cost of the postprocessing via rank-revealing information.   
 
       Therefore, in the next section, we propose a method that is iterative in nature and which has the two-fold advantage of minimizing the work for computing only what we need, in terms of a) approximating the $\bfX_{k,0}(:,j)$ where it's needed for the optimization, b) generating an approximate global basis matrix by augmenting an initial estimate with only the most non-redundant information.   The second part means we eliminate the need 
       for step 3 above in the process for computing $\bfV$, because we build $\bfV$ up to have $r$ columns as we go, rather than overbuild and and then compress to $r$ terms.          
          
%%%%%%%%%%%%%%%%%%%%%%%%%%%%%%%%%%%%%
\section{Inner-Outer Krylov Recycling}    \label{sec:recycle}

Summarizing the previous section, we see that with $\tilde{\bfA}_k := \tilde{\bfA}(\bfsfp_k)$, 
$\bfX_k := \bfX_{k,0}$, $\bfB := [\tilde{\bfB},\tilde{\bfC}]$, the relevant information for
generating the basis comes from the systems 
\begin{equation} \label{eq:blkeqn} 
%(\tilde{\bfA}(\bfsfp_k)) \bfX_{k} = [\tilde{\bfB},\tilde{\bfC}] := \bfB, 
  \tilde{\bfA}_k \bfX_k = \bfB,
\end{equation}
for several values of $k$ in order to determine the global basis.   (Note that the matrix $\bfB$ 
has a new definition from that of Section \ref{sec:background}.)  Moreover, since 
\begin{equation} \label{eq:sum}
\tilde{\bfA}(\bfsfp) = \underbrace{(\bfLL - \bfD_2 \bfG^{-1} \bfD_1)}_{\tilde{\bfA}_*} + \mbox{diag}(\mu(\bfsfp)),
\end{equation}
 the changes to the matrix as a 
function of parameter are restricted to the diagaonal.    

%The subject of this section is how to solve the sequence of systems efficiently while simultaneously and dynamically building the reduced global basis through the first steps of the optimization.      
Recycling for a sequence of systems of the form (\ref{eq:blkeqn}) was shown to be efficient in \cite{KdS2006} in
the context of optimization for shape parameters in diffuse optical tomographic imaging.   In that work, 
the authors solve (\ref{eq:blkeqn}) for all parameters selected during the course of the optimization, using different recycle spaces for each right-hand side.   
However, the focus of this paper is different since we not only want to solve the (shorter) sequence of full order problems efficiently, we need to construct this global basis, and we want  to do it in such a way that near redundant information is detected on the fly.  In this section, we first cover the basic information about recycling for a single system, then introduce our new inner-outer recycling method which avoids adding nearly redundant information to 
$\bfV$ and having to discard it (through a rank-revealing factorization) afterwards.  

%Thus, we will solve a different sequence of problems than (\ref{eq:blkeqn}) using recycling techniques in order to only update what we need.  What is particularly new in this paper is a) the problems to which the recycling is applied b) how we will build $\bfV$ from only specific parts of the solutions to the `deflated' problems.   Our approach ensures that 
%only the non-redundant information is added to the global basis matrix with each solve, so that we do not
%have the expensive proposition of computing an SVD to reduce the redundancies, as has been done previously.        

%First, we cover the basic information about recycling for a single system that is relevant to our approach.  
 
\subsection{Recycling Basics}
%\MEK{defined a macro for what we normally in recycling call ``C'', since we've already used C elsewhere.  The macro currently uses the letter K to replace C}

For this discussion, consider the linear system $\bfA \bfx=\bfb$, with symmetric
$\bfA \in \mathbb{R}^{N \times N}$ and $\bfb \in \mathbb{R}^N$. Let $\bfUU \in \mathbb{R}^{N \times n_c}$,
be given such that $\bfA \bfUU = \Co$ and $\Co^T \Co = \bfI$.  The approximate solution in
$\mrg (\bfUU)$ that minimizes the 2-norm of the residual is
\begin{equation} \label{eq:initialz} \bfz = \bfUU \bfK^T \bfb, \end{equation}
which yields the
residual $\bfr = \bfb-\Co \Co^T \bfb$ that is orthogonal to $\mrg (\Co)$.
If this solution is not adequate,
we expand the subspace as follows \cite{Stu96a}.
Let $\bfv_1 = (\bfI-\Co \Co^T)\bfb/\|(\bfI-\Co \Co^T)\bfb\|_2$ be the normalized initial residual. We use a Lanczos recurrence
with $(\bfI-\Co \Co^T)\bfA$ and $\bfv_1$ to generate the recurrence relation\footnote{We note that the matrix in the recurrence is formally equivalent to $(\bfI - \bfK \bfK^T)\bfA (\bfI - \bfK \bfK^T)$, but the right-most projector is suppressed for clarity.} 
\eqs
  (\bfI-\Co \Co^T) \bfA \bfV_m & = & \bfV_{m+1}\underbar{T}_m \Leftrightarrow \nonumber \\
  \bfA \bfV_m & = & \Co \Co^T \bfA \bfV_m + \bfV_{m+1}\underbar{T}_m ,
\label{eq:eq01}
\eqe
where $\underbar{T}_m$ is $(m+1) \times m$ tridiagonal, since $\bfA$ is symmetric.  
Next, we compute the approximate solution in $\mrg ([\bfV_m \; \bfUU])$ that minimizes
the 2-norm of the residual, $\|\bfb - \bfA( \bfV_m \bfy + \bfUU \bfz)\|_2$, as follows:
\eqs
  \min_{\bfy,\bfz} \left\| \bfb -\bfA[\bfUU \; \bfV_m ]
                         \left[\begin{array}{c}
                           \bfz \\ \bfy
                         \end{array}\right]
             \right\|_2 &  &
\nonumber \\
 = \min_{\bfy,\bfz} \left\| \bfb - [ \Co \; \bfV_{m+1}]
                         \left[\begin{array}{cc}
                          \bfI & \Co^T\bfA \bfV_m     \\
                          0 & \underbar{T}_m\\
                         \end{array}\right]
                         \left[\begin{array}{c}
                           \bfz \\ \bfy
                         \end{array}\right] \right\|_2 &  &
\nonumber \\ \label{eq:SmallMin}
=  \min_{\bfy,\bfz} \left\| \left[\begin{array}{c}
                     \Co^T \bfb \\ \xi \bfe_1
                     \end{array}\right] -
                     \left[\begin{array}{cc}
                        \bfI & \Co^T\bfA \bfV_m     \\
                        0 & \underbar{T}_m\\
                        %   0 & \underbar{H}_m
                     \end{array}\right]
                     \left[\begin{array}{c}
                           \bfz \\ \bfy
                     \end{array}\right] \right\|_2 ,
\label{eq:eqrecyc}
\eqe
where $\bfe_1$ denotes the first Cartesian basis vector in $\mathbb{R}^{m+1}$
and $\xi = \|(\bfI - \Co \Co^T) \bfb \|_2$.
The minimization in (\ref{eq:SmallMin}) corresponds to a small
least squares problem, whose solution requires only
the QR decomposition of the submatrix $\underbar{T}_m$, which
can be easily updated from the QR decomposition at step $m-1$.  The block
nature of the system means that the solution can be obtained by a three-step process:  find
$\bfy$ that solves the projected minimization problem 
\[ \min_{\bfy} \| \underbar{T}_m \bfy - \xi \bfe_1 \|_2 \]
compute $\bfz$ that satisfies (\ref{eq:SmallMin}), set $\bfx = {\bfV_m \bfy} + \bfUU \bfz $.
In reality, $\bfy_m := \bfV_m \bfy$ is computed via short term recurrences (MINRES), so the $\bfV_m$ are not
explicitly stored (see also \cite{wang2007lst, MelStu09}).   

%Ideally, the matrix $\bfU$ is such that a) it has only a relatively few columns, when storage is a premium, b) it results in a $\bfK$ such that the suitable value of $m$ required is small.   This might be the case if ..... blah blah.

\subsection{Global Basis Construction}

We start by trying to kill two birds with one stone:  if we have an estimate of the global basis matrix $\bfV$, we can investigate its use as a candidate recycling space.   If the basis can be improved, we have more full order systems to solve and we should consider solving the remaining systems by recycling.   
However, $\bfV$ may already have too many columns for it to prove computationally feasible to use as a recycle space.  

\subsubsection{Recycling on the Updated Equations}
To keep the notation simple and consistent with the previous subsection, let $k$ be fixed and set $\bfA := \tilde{\bfA}(\bfsfp_k)$, $\bfb_j = \tilde{\bfB}(:,j)$, and
$\bfx_j = \tilde{\bfX}_{k}(:,j)$.   Next, assume $\bfA \bfV = \tilde{\bfK}$, find the QR factorization 
$\tilde{\bfK} = \bfK \bfR$, and then set\footnote{In practice, $\bfUU$ is formed without inverting $\bfR$ 
explicitly.}
$\bfUU = \bfV \bfR^{-1}$, so that $\bfA \bfUU = \bfK$.  

According to the previous recycling discussion, the optimal solution in $\mrg (\bfUU) = \mrg(\bfV)$ is 
$\bfz = \bfUU \bfK^T \bfb_j = \bfV (\bfR^{-1} \bfK^T \bfb_j)$, and the initial
residual is $(\bfI - \bfK \bfK^T)\bfb_j$.  If this initial residual is small in a relative sense, then there is no need to go further, we need no iteration.   If the initial residual is not small enough, then according to the previous section,
we should expand the search space by running Lanczos to form a basis for the Krylov subspace generated by
the projected matrix $(\bfI - \bfK \bfK^T) \bfA$ and projected right-hand side $(\bfI - \bfK \bfK^T) \bfb_j$.
With this basis in the columns of $\bfV_m$, we find a an approximate solution in $\mrg([\bfUU,\bfV_m])$.  
However, unless the number of columns in $\bfV$, hence $\bfK$, is relatively small, we cannot afford to do this, because the reorthogonalization is too expensive!   Moreover, this presupposes that we actually need to find an approximation to $\bfx$ directly.  In fact, as far as updating our global basis approximation is concerned, we need only information {\it that is not already reconstructable from $\mrg(\bfV) = \mrg(\bfUU)$}.    
      
An important fact comes to light if we decompose $\bfb_j$ using the orthogonal projector $\bfK \bfK^T$: 
\begin{eqnarray} 
\bfA \bfx_j & = & (\bfI-\bfK \bfK^T)\bfb_j + \bfK \bfK^T \bfb_j  \nonumber \\  
\bfA \bfx_j - \bfK \bfK^T \bfb_j & = & (\bfI-\bfK \bfK^T)\bfb_j  \nonumber \\
\bfA \bfx_j - \bfA \bfUU \bfK^T \bfb_j & = & \bfr_j  \nonumber \\  
\bfA\underbrace{(\bfx_j - \bfUU \bfK^T \bfb_j)}_{\bfg_j} & = & \bfr_j .  \label{eq:resi} 
\end{eqnarray}
The vector $\bfg_j$ is the correction to the initial guess $\bfUU \bfK^T \bfb_j$. 
So, to obtain the incremental information to construct the global basis matrix, we should consider an iterative solution to (\ref{eq:resi}).  As we have one such system for right-hand sides $j=1,\ldots,n_{src}+n_{det}$, 
and a squence of these as we iterate over the parameters, we can employ a recycling type of approach, if we 
choose our recycle space carefully. 

Specifically, suppose $\bfr_j$ is not already suitably small, so we want to solve (\ref{eq:resi}), or alternatively, we 
want to find
\[ \min_{\bfg_j \in \mathcal{S}} \| \bfr_j - \bfA \bfg_j \|_2 ,\]
for suitable subspace $\mathcal{S}$\footnote{Note that $\mathcal{S}$ should not be $\mrg (\bfUU)$, as then the solution is zero.}.  We cannot afford to use all the columns of $\bfV$ as a recycle space in generating a $\bfV_m$ to use, because of the cost of the orthogonalization against $\mrg(K)$.  Instead, we will use a subset of the columns for right-hand side $j$ to be the recycle space,
which we'll immediately need to expand.   

Specifically, this amounts to finding $\bfUU_j \in \mathbb{R}^{N \times n_j}$, 
$\bfUU_j \subseteq\mrg (\bfUU)$, and such that $\bfA \bfUU_j = \bfK_j$, where $\bfK_j^T \bfK_j = \bfI$. 
% \begin{eqnarray*] (\bfI - \bfK_j \bfK_j^T) \bfA \bfx_ & = & (\bfI - \bfK_j \bfK_j)^T \bfr_j \\
 %                                                  & = & (\bfI - \bfK \bfK^T) \bfb_j  \end{eqnarray*}
Now we use $\mathcal{S} =\mrg([\bfUU_j,\bfV_m^{(j)}])$ where the $\bfV_m^{(j)}$ are the Lanczos vectors for $K_m\left( (\bfI - \bfK_j \bfK_j^T) \bfA, (\bfI - \bfK_j \bfK_j^T)\bfr_j \right)$ (compare to (\ref{eq:eq01})).%  to look for a solution.  
That is, we want to solve
\[ \min_{\bfz,\bfy} \left\| \bfr_j - \bfA [\bfV_{m}^{(j)}, \bfUU_j] \bea{c} \bfy \\ \bfz \eea  \right\|.\]
 %Comparing to (\ref{eq:eq01}), we will generate a recurrence for the projected (with projector $\bfI - \bfK_j \bfK_j^T)$) problem.  
 Importantly, this choice for $\bfUU_j$ gives $(\bfI - \bfK_j \bfK_j^T) \bfr_j = \bfr_j$.  
 % that if we apply the projector $(\bfI - \bfK_j \bfK_j^T)$ to both sides of (\ref{eq:resi}), the right-hand side
 %$\bfr_j$ remains unchanged  because of the relationship between $\bfK_j$ and $\bfK$.  
 Thus, we observe $\bfv_1^{(j)} = \bfr_j/\| \bfr_j \|$.  Next we use the Lanczos recurrence
with $(\bfI-\bfK_j \bfK_j^T)\bfA$ and $\bfv_1^{(j)}$ to generate the recurrence relation
\eqs
  (\bfI-\bfK_j \bfK_j^T) \bfA \bfV_m^{(j)} & = & \bfV_{m+1}^{(j)} \underbar{T}_m^{(j)} \Leftrightarrow \nonumber \\
  \bfA \bfV_m^{(j)} & = & \bfK_j \bfK_j^T \bfA \bfV_m^{(j)} + \bfV_{m+1}^{(j)} \underbar{T}_m^{(j)} .
\eqe
Now $\bfy, \bfz$ are found by (compare to (\ref{eq:eqrecyc})) solving
\[ \min_{\bfy,\bfz} \left\| \bea{c} 0 \\ \xi \bfe_1 \eea -
                                         \bea{cc} \bfI & \bfK_j^T \bfA \bfV_m^{(j)} \\
                                                         0 & \underbar{T}_m^{(j)} \eea \bea{c} \bfz \\ \bfy \eea \right\|_2 . \] 

Hence, $\bfg_j = {\bfV_m^{(j)} \bfy} + \bfUU_j \bfz$, where 
$\bfz = -\bfK_j^T \bfA \bfV_m^{(j)} \bfy$ and $\bfy_m^{(j)}:= \bfV_m^{(j)} \bfy$ is generated by a short term recurrence.  So $\bfx_j = \bfy_m^{(j)} - \bfUU_j \bfK_j^T \bfA \bfy_m^{(j)} + \bfUU \bfK^T\bfb_j$. Since $\mrg(\bfUU_j) \subseteq \mrg(\bfUU) = \mrg(\bfV)$, the relevant incremental information about $\bfx_j$ that cannot already be expressed using the columns of $\bfV$ is $\bfV_m^{(j)} \bfy$.   Therefore, we extend the global basis with this vector.  We repeat this process for any $j$ for which the initial residual $\bfr_j$ is not already small enough\footnote{The integer $m$ for which the solution estimate is good enough will 
  vary depending on the system -- that is, $m = m_j$ -- but for ease in notation we have omitted the subscript on $m$.}.   At most, for a fixed $k$, we have added $n_{src} + n_{det}$ columns to 
  $\bfV$, but in theory, we may add substantially fewer.        
                                                   
%Our goal is to solve 
%\[ \min_{\bfg_j \in \mrg([\bfV_m,\bfU_j]) } \| \bfr_j - \bfA \bfg_j \|_2 .\]
%The preceding subsection shows that this is obtained by first computing $\bfy$, the solution to the projected problem, then we would have $\bfg_j = \underbrace{\bfV_m y}_{\bfy_m} + \bfU_j \bfz$. 

%We should note that although we are only interested in $\bfg_j$ for augmenting $\bfV$, it is possible to back out
%$\bfx_j$.  
%From ~\ref{eq:eqrecyc}, we have that $\bfz = \bfK_j^T \bfr - \bfK_j^T \bfA \bfy_m$, which means, 
%\begin{eqnarray*}
%\bfg = \bfx - \bfV \bfK^T \bfb & \approx & \bfUU_j \bfz + \bfy_m \\
%\bfx - \bfV \bfK^T \bfb & \approx & \bfUU_j \bfK_j^T \bfr - \bfUU_j \bfK_j^T \bfA \bfy_m + \bfy_m \\
%\bfx & \approx & \bfV \bfK^T \bfb - \bfUU_j \bfK_j^T \bfA \bfy_m + \bfy_m. \\
%\end{eqnarray*}

Now we consider the remaining issues: a) what to use as the initial guess to the global basis $\bfV$ and b) how to choose the individual columns to specify the $\bfUU_j$ when we do need to solve (\ref{eq:resi}).      

\subsubsection{Identifying the Recycling Spaces}
In \cite{KdS2006} when working on 
a different parametric inverse model problem for DOT, the authors observed that a good recycle subspace for
one right-hand side did not necessarily make a good recycle space for the next right-hand side.  Instead, they employed a different recycle space for each right-hand side in the sequence, with the caveat that the recycle spaces did have a common subspace pertaining to a specific (approximate) invariant subspace.   That invariant subspace 
corresponded to the smallest eigenvalues, because that subspace remained relatively unchanged through the optimization process.

Even though in the present paper we are using a different image parameterization than what was used
in \cite{KdS2006}, we also observe that the invariant subspace due to the smallest several eigenvalues
for the $\tilde{\bfA}(\bfsfp_k)$ remains unchanged.  Thus, we adopt the same approach here. We use a different recycle space $\bfUU_j$ for each right-hand side, but seed each with the same invariant subspace, plus right-hand-side specific information.

First, we compute (approximate) eigenvectors of $\tilde{\bfA}_0$ that correspond to the smallest eigenvalues.  The small eigenvalues of the $\tilde{\bfA}_k$ matrices remain close from one system to the next suggesting that the corresponding invariant subspaces also remain close.  We refer the reader to \cite{KdS2006} for more details and theory. We have found experimentally that $10$ eigenvectors is sufficient for the invariant subspace, while keeping the recycle space small.  So, we set $\bfUU_0 \in \mathbb{R}^{n \times 10}$ to contain (estimates of)
those 10 vectors.  Spending the computational resources to get an accurate invariant subspace so that it can
be deflated from the right-hand side has been observed to be worthwhile in other large-scale applications, such as QCD, as well \cite{Stath2007,Stath2009}.

%\MEK{Eric, you wanted us to add ref to Waisman, Fish, et al, 2004, but where?}
Initially, we let $\bfV = ([\bfUU_0, \bfX_0])$ and $\bfUU_j = [\bfUU_0,\bfX_0(:,j)]$.   Note that   
$\bfV$ has the initial solutions to all the right-hand sides, while $\bfUU_j$ has only the solution from the $j^{\rm th}$ right-hand side.  If we find that $\bfr_j$ is not suitably small, perform the recycling outlined above on (\ref{eq:resi}), and we update both $\bfV$ and $\bfUU_j$ with $\bfy_m^{(j)}$.  The $\bfy_m^{(j)}$ is appended to $\bfV$ every time we need to do recycling, while it is only appended to $\bfUU_j$ if we are working on the $j^{\rm th}$ right-hand side.  This ensures that $\bfV$ contains information pertinent to the entire system, while $\bfUU_j$ is kept small.

\subsection{The Algorithm} \label{subsec:alg}
Algorithm \ref{alg:recycbasis} describes our dynamic process.  Details on efficient
implementation of various steps in the algorithm will be addressed in the next section.    

\begin{algorithm}[ht]
\SetAlgoLined 
$\bfUU_0 \Leftarrow 10 \: \mbox{eigenvectors of} \: \tilde{\bfA}_0$, $\bfX_0$ solves $\tilde{\bfA}_0 \bfX_0 = \bfB$ \\
$\bfV \Leftarrow \text{basis for }\mrg([ \bfUU_0, \bfX_0])$ \\
$\bfUU_j \Leftarrow [\bfUU_0, \bfX_0(:,j)]$ \\
\For{$i = 1:K$ \% for each interpolation point $i$}{ 
	\For{$j = 1:nrhs$}{
      	\% Check if $\bfV$ is a good enough space \\
      	$\tilde{\bfK} = \tilde{\bfA}_i \bfV$ \\ \label{KAV}
      	$[\bfK,\bfR] = qr(\tilde{\bfK},0)$ \\	\label{QRV}
      	$\bfV = \bfV / \bfR$ \%implicit only; now $\bfUU,\bfV$ same \\	\label{VR}
      	$\bfr_j = \bfB(:,j) - \bfK \bfK^T \bfB(:,j)$ \\
          \If{$\frac{||\bfr_j||}{||\bfB(:,j)||} > tol$}{
          	\% MINRES recycling using $\bfUU_j$ \\
          	$\tilde{\bfK}_j = \tilde{\bfA}_i \bfUU_j$ \% have already done this product\\	\label{KAU}
          	$[\bfK_j,\bfR] = qr(\tilde{\bfK}_j,0)$ \% need not be done from scratch\\	\label{QRU}
          	$\bfUU_j = \bfUU_j / \bfR$ \\		\label{UR}
          	Solve $(\bfI - \bfK_j \bfK_j^T)\tilde{\bfA}_i\bfy_m^{(j)} = \bfr_j$ with MINRES \\
            $\bfV \Leftarrow [\bfV, \bfy_m^{(j)}]$ \\
            $\bfUU_j \Leftarrow [\bfUU_j, \bfy_m^{(j)}]$ \\
          }
	}
}  
\caption{Recycling and Global Basis Construction\label{alg:recycbasis}}    		
\end{algorithm}

\section{Algorithm Analysis}    \label{sec:analysis}

Algorithm \ref{alg:recycbasis} given in the previous section describes how we solve the sequence of systems while generating the approximate global basis.  Our approach not only solves the sequence of systems efficiently but also builds the reduced global basis with only non-redundant information, therefore eliminating unnecessary solves as well as the need for an expensive SVD and a corresponding ad hoc approach for truncation.  Additionally, we can keep the cost down by exploiting the fact that we update $\bfV$ and $\bfUU_j$ one column at a time.

Next, we explain how our approach differs from \cite{Gugetal2015}.  We also show that it is different from simply applying the recycling as in \cite{KdS2006} to solve the systems and then doing an SVD to get the global basis.

\subsection{System Solves}

Our new approach is an improvement over Algorithm \ref{alg:alg1}.  To see this, we consider implementation of Step 1.
%As mentioned, the approach for generating the global basis that we want to improve upon is the approach in \cite{Gugetal2015}, which suggests solving (\ref{eq:blkeqn}) for several choices of $\bfsfp_i$, concatinating the
%solutions $\bfX_k$, and compressing by truncating an SVD.  To implement this and 
To find the $\bfX_k$, one can of course use MINRES directly for each right-hand side.  Alternatively, one can use recycling with the $\bfUU_j$ for the respective right-hand side, across all the systems\footnote{This is essentially the approach in \cite{KdS2006}, with recycle spaces having common invariant subspace information but tailored to the particular right-hand side.   But they do further tuning of the recycle spaces to account for where one is in the optimization process.}.   
The recycling would consist of generating $\hat{\bfV}_m^{(j)}$ 
as a basis for ${\mathcal K}_m( (\bfI - \bfK_j \bfK_j^T) \bfA, (\bfI - \bfK_j \bfK_j^T) \bfb_j)$, with the intent of approximating $\bfx_j$ over $\mrg([\bfUU_j,\hat{\bfV}_m^{(j)}])$.   

In contrast, in our new approach, we solve (\ref{eq:resi}).  Now $\bfV_m^{(j)}$ is generated as a basis
for ${\mathcal K}_m( (\bfI - \bfK_j \bfK_j^T)\bfA, (\bfI - \bfK \bfK^T)\bfb_j)$.   The two Krylov spaces differ in the
two approaches by the right-hand sides.        
Also,  
we want to approximate $\bfx_j - \bfx_{0,j}$, with $\bfx_{0,j} = \bfUU \bfK^T \bfb_j$ rather than $\bfx_j$.  
A further numerical comparison is provided in Section (\ref{subsec:numcomp})
 
It should be clear that recycling using $\bfUU_j$ must have some advantage over not recycling at all.   In the first
place, since $\bfUU_j$ contains an approximate invariant subspace, MINRES convergence on the projected systems would behave as if part of the spectrum has been deflated.

Additionally, the right-hand sides of (\ref{eq:resi}) are residuals that have been already made small across spectral components {\it other than just those included in the invariant subspace}, since these residuals are $\bfb_j$'s orthogonalized against the entire $\bfK$, not just the (much smaller) $\bfK_j$.    
An argument for why $\bfr_j$ is small in norm follows along the lines of subsection $5.2$ in \cite{KdS2006}, and assumes 
that $\tilde{\bfA}_0 - \tilde{\bfA}_1$ is small over the invariant subspace of $\tilde{\bfA}_0$ corresponding to the
smallest eigenvalues (smooth modes), and makes use of the fact that the columns of $\bfX_0$ are relatively smooth.  
%We will investigate for $k=1$, for the first right-hand side.  
If we assume that the values of absorption in our object and in the background are known and optimize only for the shape parameters, $\tilde{\bfA}_0 - \tilde{\bfA}_1$ is both diagonal, possibly low rank, with smooth modes made small by the operator.  %As the argument requires only minor adjustment for our use of multiple columns of $\bfX_0$, and is therefore omitted.   

We claim that using space $\bfUU_j$ instead of the 
(expensive) $\bfV$ does not require many more iterations.  
%Due to the cost, use of $\bfUU_j$ is therefore preferable as seen in Section ~\ref{subsec:numcomp}.  
%Note that the $\bfUU_j$ used in \cite{KdS2006} and the $\bfUU_j$ used in our approach differ following after the first update (see also the numerical results comparison in Section \ref{subsec:numcomp}) .  
%Here, we use a different right-hand-side and therefore, although both $\bfUU_j$ are updated with a solution to the projected problem, these solutions will differ and therefore the $\bfUU_j$ will differ after the first update.     

\subsection{Global Basis}
\begin{figure}[h]
	\centering
	\includegraphics[scale=0.75]{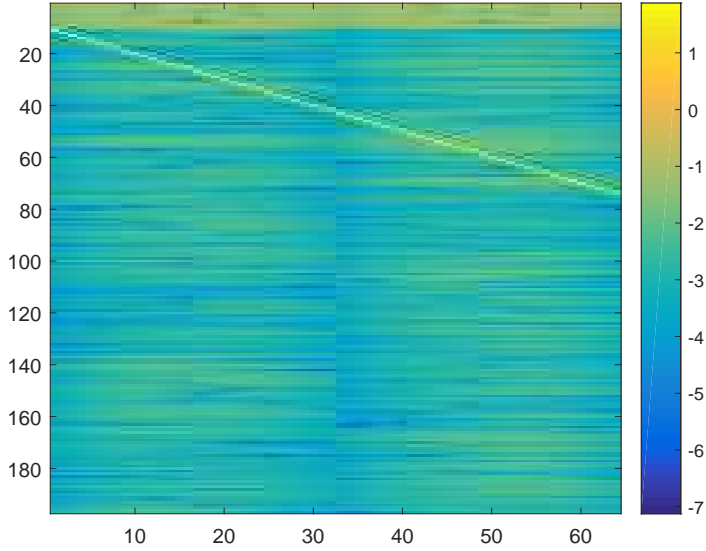}
	\caption{Logarithmically scaled image of the absolute values of the coefficients of solutions $\bfX_{5}$ in the directions of $\bfV$, where $\bfV$ has 197 columns.  Recall that $\bfX_5$ has 32 columns to
	source positions, and 32 columns corresponding to receiver positions.  The first 10 columns of $\bfV$ contain
	the invariant subspace of $\tilde{\bfA}_0$ corresponding to its smallest 10 eigenvalues, the next 64 columns 
	correspond to $\bfX_{0}$, and the remaining columns have been constructed using the update procedure
	in Algorithm \ref{alg:recycbasis} for $K_*=3$.} 
	\label{fig:amXComponents}
\end{figure}

%The second goal of the algorithm is to generate the reduced global basis.  
It is natural to ask  why $\bfV$ is a good global basis for the ROM problem.  
In addition to solution information, our $\bfV$ contains an (approximate) invariant subspace corresponding to the smoothest modes for $\tilde{\bfA}_0$.  % But $\bfX_k$ are expected to be smooth.  
By Theorem $4.1$ in \cite{KdS2006}, we know that if the changes to $\tilde{\bfA}_k$ are
concentrated over the high frequency modes, the invariant subspace consisting of eigenvectors corresponding to the smallest eigenvalues (low frequency modes) remain close.   The $\bfX_k$ are expected to be smooth.     
From a ROM standpoint, it would only be helpful to include this information in $\bfV$ if we expect the solutions to be well represented in terms of the smooth modes.  
%Since $\tilde{\bfA}^{-1}$ acts as a blurring operator and the smallest eigenvalues represent smooth modes, we expect this to be the case.  
Figure ~\ref{fig:amXComponents} shows, in log scale, the absolute values of the coefficients of the solutions $\bfX_5$ in the directions of $\bfV$ for Experiment $1$ in Section~\ref{sec:numerical}.  Note that these solutions were not used to build $\bfV$.  You can see that the solutions have large components in the directions of the invariant subspace of $\tilde{\bfA}_0$, as well as the corresponding column of $\bfX_0$.  

\subsection{Implementation Issues} \label{subsec:cost}

In subsection~\ref{subsec:alg}, we gave an algorithm for constructing a global basis using the full order model solves using recycling. We will now discuss the cost of Algorithm \ref{alg:recycbasis} and how adding one column at a time to $\bfV$ and $\bfUU_j$ helps to keep the cost down.

We will assume that $\bfUU_0$ is a matrix of $k_u$ columns containing a basis for the invariant subspace of $\tilde{\bfA}_0$ corresponding to the smallest eigenvalues and that $\bfX_0$ is known.  Note that $\bfX_0$, $\bfUU_0$ can be precomputed off-line and can be reused for
other experiments.  %There is a cost for solving $\tilde{\bfA}_0 \bfX_0 = \bfB$ and for generating the approximate invariant subspace, $\bfUU_0$, but for a class of experiments, if $\tilde{\bfA}_0$ is always the initial matrix, $\bfX_0$ and $\bfUU_0$ can be precomputed off-line.  
Again, we will use the 
fact (see (\ref{eq:sum})) that $\tilde{\bfA}_k = \tilde{\bfA}_* + \bDelta_k$, where the first term in the sum is fixed and
$\bDelta_k = \mbox{diag}(\mu(\bfsfp_k))$.   

Let's consider solving system $1$ using Algorithm~\ref{alg:recycbasis}.  In Algorithm~\ref{alg:recycbasis} line~\ref{KAV}, $\tilde{\bfK} \coloneqq (\tilde{\bfA}_* + \bDelta_1)[\bfUU_0, \bfX_0]$. We can precompute and save the products $\tilde{\bfA}_* \bfX_0$ and $\tilde{\bfA}_* \bfUU_0$, so only updating by $\bDelta_1 \bfUU_0$ and $\bDelta_1 \bfX_0$ is required.  

Next, in line~\ref{QRV} we need to compute the QR factorization of $\tilde{\bfK}$. We will compute this in a way that does not require full re-orthogonalization each time it is repeated. First, partition $\tilde{\bfK} = [\tilde{\bfK}_a, \tilde{\bfK}_b]$, where the first block corresponds to the number of columns of $\bfUU_0$. Next compute, $\bfQ_1 \bfR_1 = \tilde{\bfK}_a$.  Then, compute $(\bfI - \bfQ_1 \bfQ_1^T)\tilde{\bfK}_b = \tilde{\bfK}_b - \bfQ_1(\bfQ_1^T \tilde{\bfK}_b) = \bfQ_2 \bfR_2$. So, $\tilde{\bfK}_b = \bfQ_2 \bfR_2 + \bfQ_1(\bfQ_1^T \tilde{\bfK}_b)$. It follows that a QR factorization is \begin{align*}
   [\tilde{\bfK}_a,\tilde{\bfK}_b] &= [\bfQ_1 \bfR_1, \bfQ_1 (\bfQ_1^T \tilde{\bfK}_b) + \bfQ_2 \bfR_2 ] \\
                                   &= [\bfQ_1,\bfQ_2] \left[ \begin{array}{cc} \bfR_1 & \bfQ_1^T \tilde{\bfK}_b \\ 0 & \bfR_2 \end{array} \right] .
\end{align*} 
We call $\bfV = [\bfUU_0, \bfX_0]\bfR^{-1}$, noting that $\bfR^{-1}$ need not be applied explicitly.  Now we have $(\tilde{\bfA}_* + \bDelta_1)\bfV = \bfK$, and $\bfK$ has orthonormal columns.  

Suppose we do recycling for system $1$ and right-hand side $1$. In Algorithm~\ref{alg:recycbasis} line~\ref{KAU}, we need $\tilde{\bfK}_1 \coloneqq (\tilde{\bfA}_* + \bDelta_1)[\bfUU_0, \bfX_0(:,1)]$.  We have already formed this product above, so we just have to select the right columns of $\tilde{\bfK}$. Moreover, the QR factorization of $\tilde{\bfK}_1$ is computed from $\bfQ_1$ and $\bfR_1$. All we need to compute is $\tilde{\bfK}_1(:,k_u+1) - \bfQ_1(\bfQ_1^T\tilde{\bfK}_1(:,k_u +1))$ and then normalize it.  The normalization constant becomes the lower right corner component of the upper triangular matrix.  

Next, we solve the projected problem with MINRES and append the solution, $\bfy_m^{(1)}$, to $\bfV$ and $\bfUU_j$. We will check to see if the newly enlarged $\bfV$ is sufficient to represent the solution for the second right-hand side. In order to do this, we compute 
\[(\tilde{\bfA}_* + \bDelta_1)[\bfV, \bfy_m] = [(\tilde{\bfA}_* + \bDelta_1)\bfV, (\tilde{\bfA}_* + \bDelta_1)\bfy_m] = [\bfK, \mathbf{z}].\]
Since $\bfK$ already has orthogonal columns, we need only to compute 
%finding the QR decomposition of $[\bfK, \mathbf{z}]$ is trivial. Using 
%$\mathbf{z} - \bfK(\bfK^T \mathbf{z}) = \rho \bfq$,  
\[ [\bfK,\mathbf{z}] = [\bfK,\bfq] \underbrace{\left[ \begin{array}{cc} \bfI & \bfK^T \mathbf{z} \\ 0 & \rho \end{array} \right]}_{\hat{\bfR}} ,\]
where $\mathbf{z} - \bfK(\bfK^T \mathbf{z}) = \rho \bfq$,
so we have $\bfK \leftarrow [\bfK,\bfq], \bfV \leftarrow [\bfV, \bfy_m] \hat{R}^{-1}$. 

For a recycle solve for right-hand side 2, we follow the same procedure.  For additional right-hand sides, only incremental new calculations are needed.

\section{Numerical Results}  \label{sec:numerical}
We present two experiments on a $201 \times 201$ mesh (which gives us $40401$ degrees of freedom) for the forward problem.  We use $32$ sources and $32$ detectors in the model, meaning that $\bf{B}$ in (\ref{eq:blkeqn}) will have 64 columns.  The image space is parameterized using parametric level sets (PaLS) (see \cite{Aghasi_etal11} for details).  We use $25$ compactly supported radial basis functions to define the PaLS image, which results in a total of 100 parameters for the optimization problem (\ref{eq:InvProb}). 

  The ground truth images for Experiments 1 and 2 are given in Figure \ref{fig:amTrue} and \ref{fig:cupTrue}, respectively.  We note that these images cannot be exactly reconstructed via the image space parameterization
  we are using. Thus,
we avoid the so-called inverse crime.   To obtain the noisy data we added 1\%, noise to the simulated true measured data in each of our experiments.  

\begin{figure}
	\centering
	\includegraphics[scale=0.3]{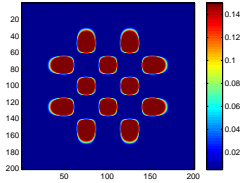}
	\caption{Initial Absorption Image}
	\label{fig:initialAbs}
\end{figure}

We
solve
the optimization problems using the TREGS \cite{StuKil11c} algorithm.  We stop the optimization when the residual norm falls below 1.1 times the noise level.  We report and compare the results
for two cases in each experiment.  First, we report results assuming that the full order problem was used to compute the function and Jacobian evaluation at each step.   Then we report results using the ROM to replace the function and Jacobian evaluation.  Figure~\ref{fig:initialAbs} gives the absorption image using the initial set of parameters.  We used $K_* = 3$ systems for each experiment to create the reduced order model space.  
The tolerance in line 11 of Algorithm 1 was set to be $10^{-7}$.  All of the experiments were run using a laptop with a $3.20$ GHz processor and $16.0$ GB RAM using MATLAB R2014a.  

\subsection{Experiment 1}
In Experiment $1$, we needed to solve $187$ large, single right-hand side systems to generate what we needed to construct the global basis matrix (note that the 64 of these corresponding to $\bfX_0$ could have been pre-computed off-line).  Including the additional $10$ eigenvectors of $\tilde{\bfA}_0$ that were used as the first 10 columns of $\bfV$,  $\bfV$ has 197 columns and thus the reduced model has order $197$.  Therefore, the reduced models require solutions to linear systems of size $197 \times 197$ rather than $40401 \times 40401$ for the full order model.  

The optimization using the full order model for Experiment $1$ required $30$ function evaluations and $15$ Jacobian evaluations.  In comparison, the optimization run using the reduced order model, once it's been generated, for function and Jacobian evaluations required $28$ function evaluations and $14$ Jacobian evaluations, indicating that using a ROM in place of FOM does not greatly impact convergence rate of the optimization.  
%The work summary is given in Table \ref{tab:worksummary}.   
The bottom line is that
solving the optimization using the full order model requires the solution of 1440 systems of size $40401 \times 40401$.  On the other hand, solving using our approach requires solution of 187 systems of size $40401 \times 40401$, which are used to construct $\bfV$ during the first few optimization steps.  The remainder of the work is in solving systems of size $197 \times 197$ until the convergence tolerance for the optimization is achieved.

%\begin{table}
%\begin{center}
%\begin{tabular}{|l||c|c|}
% \hline  & FOM  &  ROM \\  \hline
%system size &  40401 & 197 \\
%function evals & 30 & 28 \\
%Jacobian evals & 15  & 14 \\ 
%total no. FOM single-RHS system solves & 1440 &  187 \\  \hline
%total no. ROM single-RHS system solves & -- & 42 \\ \hline 
%\end{tabular}
%\caption{\label{tab:worksummary}Summary of work comparison, using ROM vs. using FOM in the solution of the optimization problem.   This does not include the time to compute/estimate the 10 eigenvectors of $\tilde{\bfA}_0$
%used to generate the ROM, a computation which can be done off-line and reused if other DOT images are to be recovered.}
%\end{center}
%\end{table}    

\begin{figure}[h]
        \centering
        \begin{subfigure}[t]{.3\textwidth}
        		\centering
                \includegraphics[width=4cm]{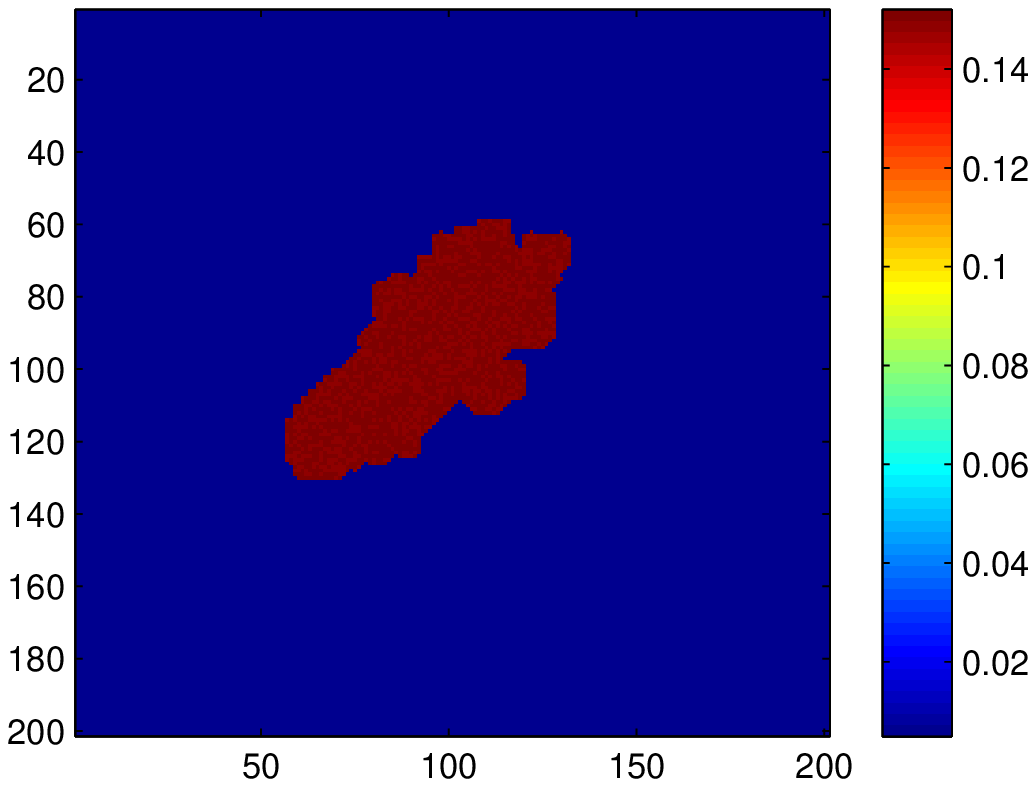}
                \caption{Original anomaly with normally distributed noise added.}
                \label{fig:amTrue}
        \end{subfigure}
        \begin{subfigure}[t]{.3\textwidth}
        		\centering
                \includegraphics[width=4cm]{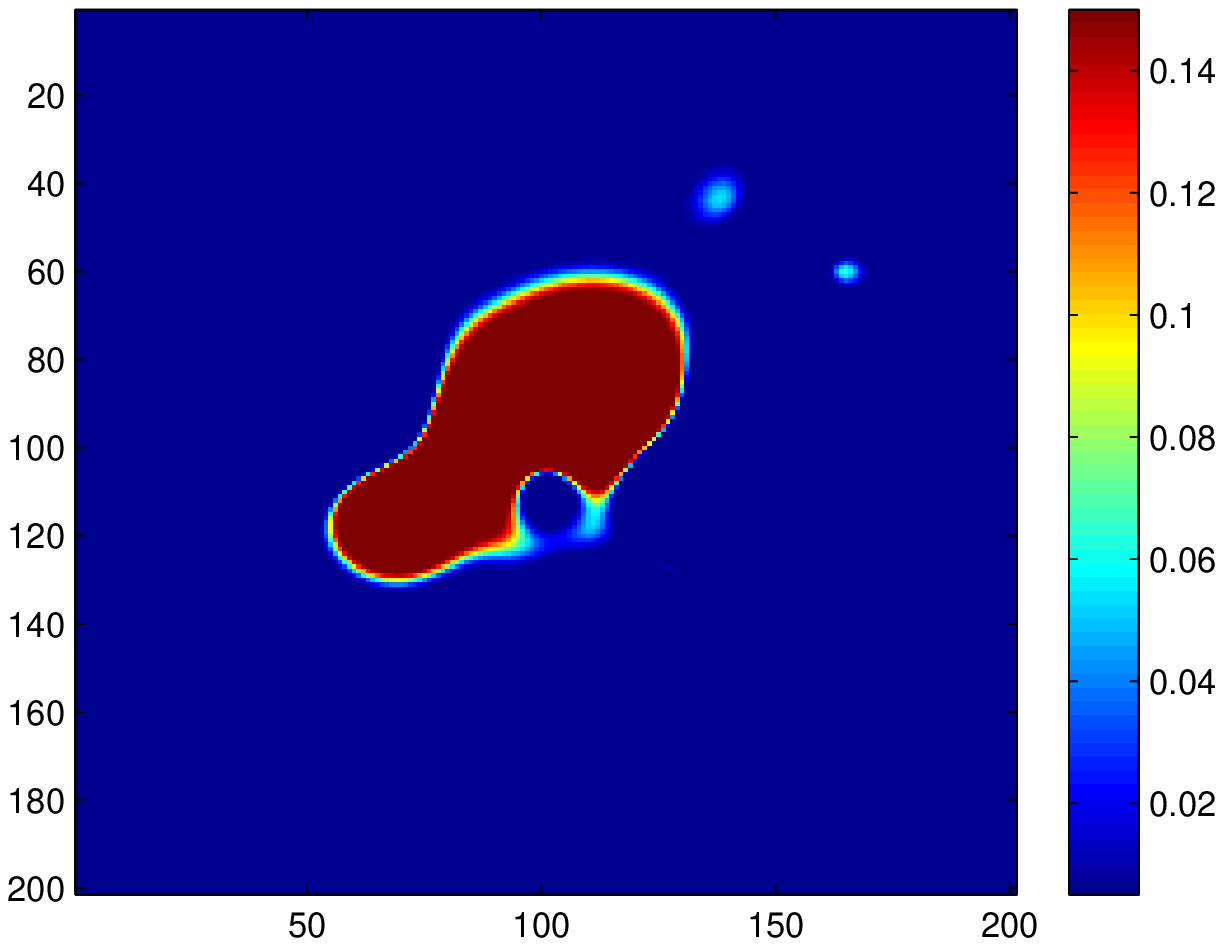}
                \caption{Reconstruction using the full order model.}
                \label{fig:amFOM}
        \end{subfigure}
        \begin{subfigure}[t]{.3\textwidth}
        		\centering
                \includegraphics[width=4cm]{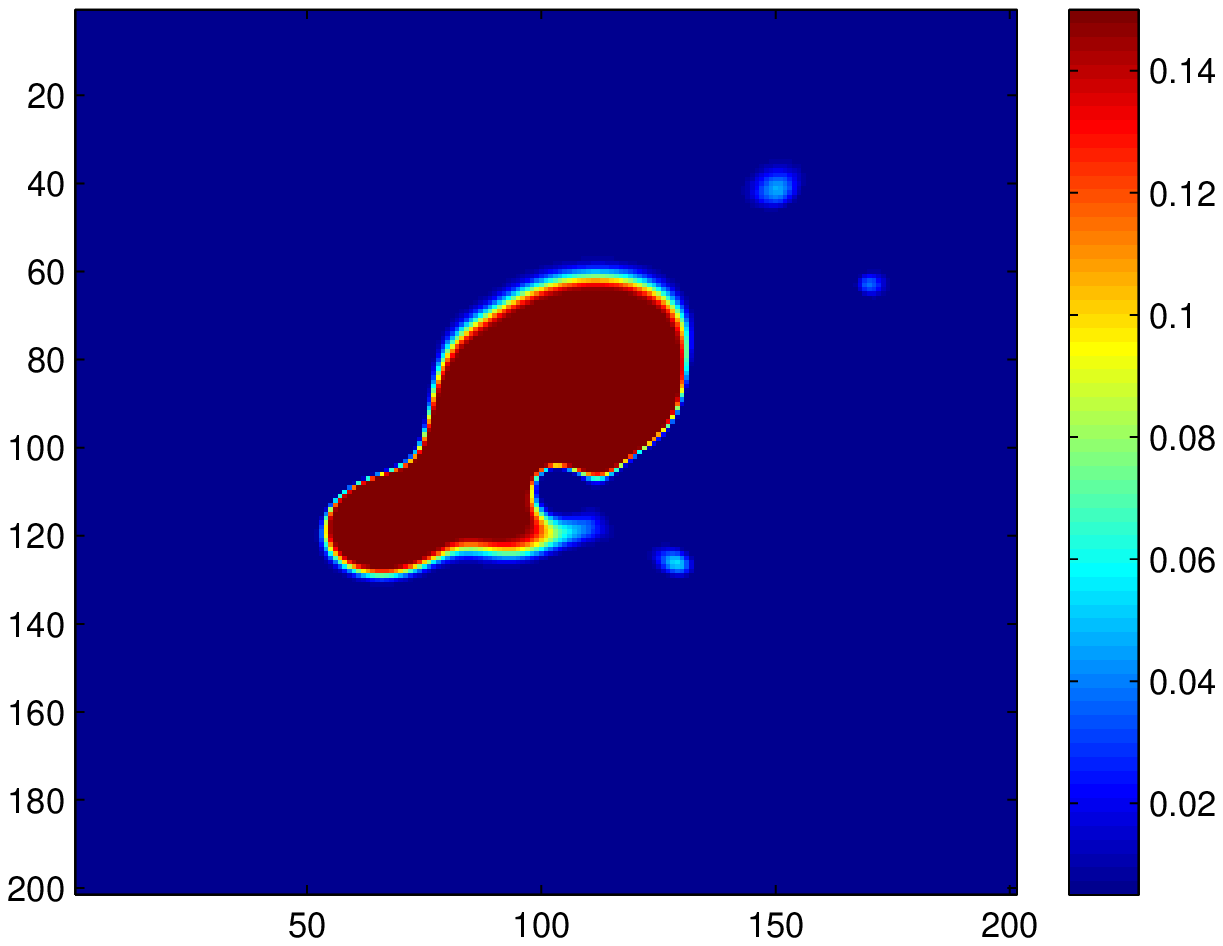}
                \caption{Reconstruction using the reduced order model.}
                \label{fig:amROM}
        \end{subfigure}
        \caption{Results for Experiment 1.  Reconstruction on a $201 \times 201$ mesh, resulting in $40401$ degrees of freedom in the forward model and $197$ degrees of freedom in the reduced model for the forward model.  $32$ sources, $32$ detectors, and $25$ basis functions were used. }\label{fig:am}
\end{figure}

Figure~\ref{fig:am} shows the reconstructions for Experiment $1$.  Figure~\ref{fig:itsTable} also includes the number of (unpreconditioned) MINRES iterations for each experiment with and without recycling.  Although the tables only show a sample of results, it is clear that the iterations decrease from one right-hand side to the next, and system to system, using our approach.  The jump in number of iterations for right-hand-side $33$ comes from the fact that we concatenated $\tilde{\bfB}$ and $\tilde{\bfC}$ to form one right-hand-side for the symmetric transfer function, so the 33rd right-hand side corresponds to the first column in $\tilde{\bfC}$.     

\begin{table}[h]
        \centering
				\begin{tabular}{| c | c | c | c | c | c |}
				    \cline{1-6}
				      &  & \multicolumn{2} { c| } {Experiment 1} & \multicolumn{2} { c| } {Experiment 2} \\ \hline
				      \multirow{2}{*}{}
				    System & RHS & MINRES Its & MINRES Its & MINRES Its & MINRES Its  \\ 
				     &  &  & with Recycling &  & with Recycling  \\
				    \hline
				    \multirow{6}{*}{1} & 1 & 463 & 140 & 470 & 127\\
				    
				    & 20 & 541 & 52 & 506 & 53 \\  
				    & 32 & 487 & 0 & 493 & 0 \\ 
				    & 33 & 467 & 124 & 470 & 120 \\ 
				    & 53 & 528 & 57 & 514 & 48 \\ 
				    & 64 & 489 & 0 & 494 & 0 \\ \hline
				    \multirow{6}{*}{2} & 1 & 474 & 118 & 501 & 124 \\ 
				    & 20 & 513 & 35 & 545 & 38 \\ 
				    & 32 & 497 & 0 & 526 & 5 \\ 
				    & 33 & 474 & 105 & 500 & 132 \\ 
				    & 53 & 526 & 37 & 567 & 127 \\  
				    & 64 & 497 & 0 & 532 & 0 \\ \hline
				    
				    \end{tabular}
        \caption{Number of MINRES iterations for Experiments 1 and 2.}
        \label{fig:itsTable}
\end{table}

\subsection{Experiment 2}
In Experiment $2$, a total $188$, 40401 x 40401 single-right-hand side systems were solved to compute our global basis.  The reduced order model has order $198$.  Therefore, the reduced models require solutions to linear systems of size $198 \times 198$ rather than $40401 \times 40401$ for the full order model.  The optimization using the full order model required $126$ function evaluations and $78$ Jacobian evaluations.   The optimization run using our reduced order model took $123$ function evaluations and $76$ Jacobian evaluations to converge to our stopping criterion,  so again, there is no negative impact on convergence rate by replacing the FOM with the ROM. 
The difference in the total number of large (40401 x 40401) single-right-hand side systems that need to be solved, though, is even more pronounced in
this example than in the last:  6,528 are needed for the FOM approach vs. only 188 for the ROM approach.  Moreover, the work involved in solving for the latter systems is reduced, since MINRES requires fewer iterations
due to the recycling.   

\begin{figure}[h]
        \centering
        \begin{subfigure}[t]{.3\textwidth}
        		\centering
                \includegraphics[width=4cm]{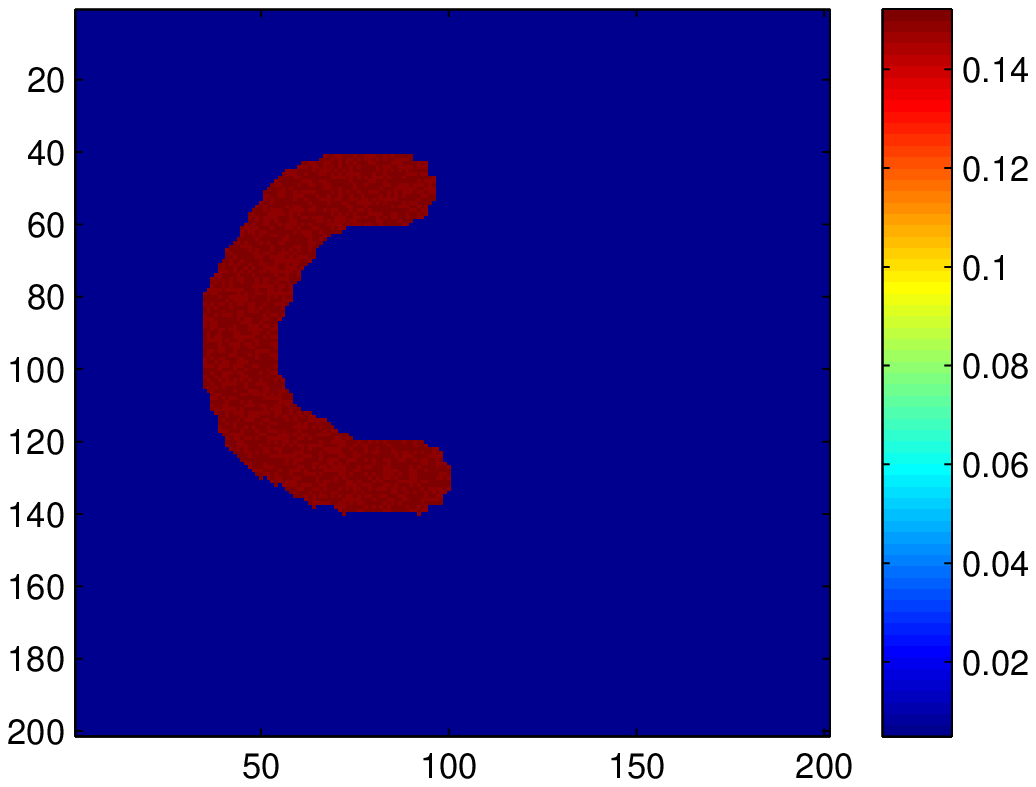}
                \caption{Ground truth image.}
                \label{fig:cupTrue}
        \end{subfigure}
        \begin{subfigure}[t]{.3\textwidth}
        		\centering
                \includegraphics[width=4cm]{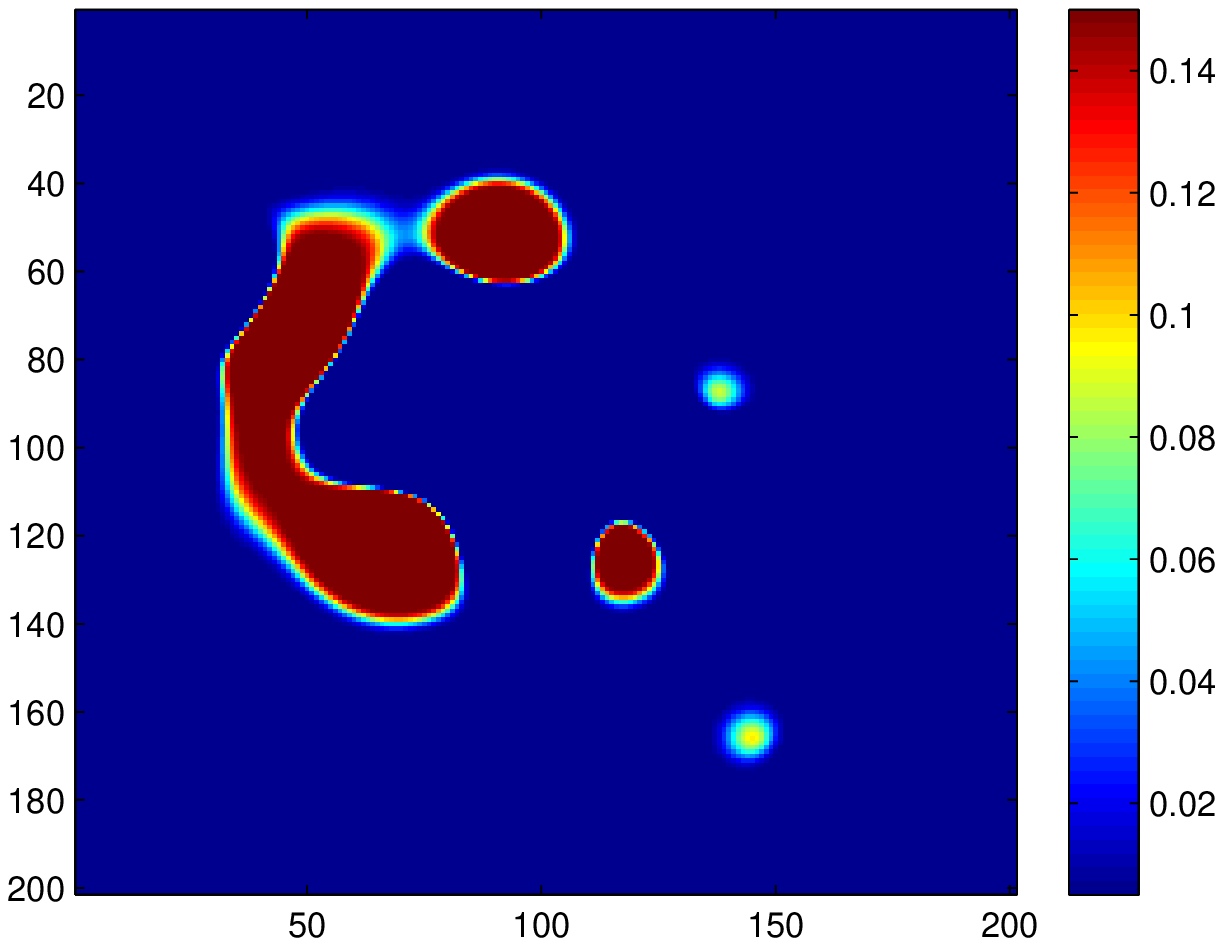}
                \caption{Reconstruction via the FOM.}
                \label{fig:cupFOM}
        \end{subfigure}
        \begin{subfigure}[t]{.3\textwidth}
        		\centering
                \includegraphics[width=4cm]{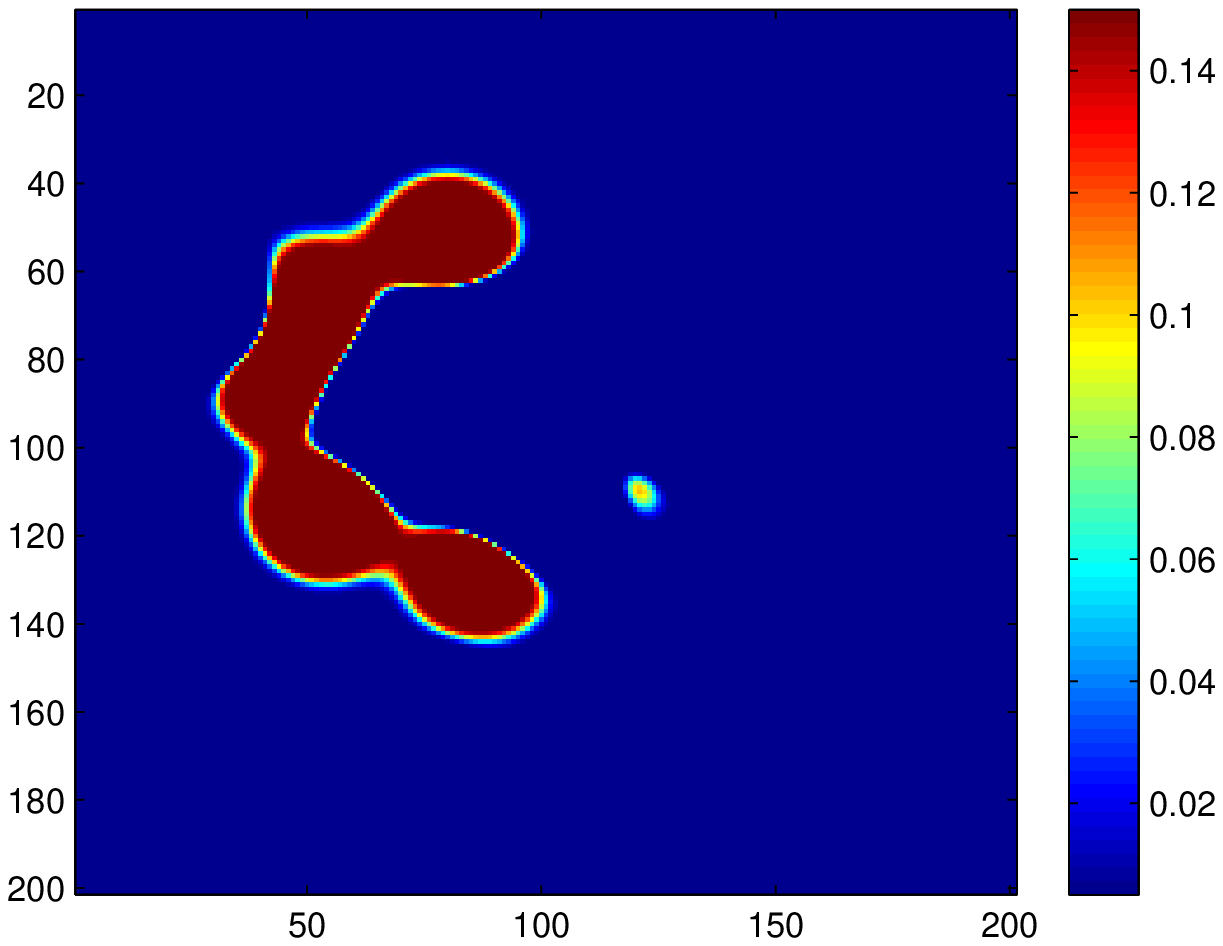}
                \caption{Reconstruction via the ROM.}
                \label{fig:cupROM}
        \end{subfigure}
        \caption{Results for Experiment 2.  Reconstruction on a $201 \times 201$ mesh, resulting in $40401$ degrees of freedom in the forward model and $198$ degrees of freedom in the reduced order model.  $32$ sources, $32$ detectors, and $25$ basis functions were used. }\label{fig:cup}
\end{figure}

  Figure~\ref{fig:cup} shows the reconstructions for Experiment $2$.  Again, Figure~\ref{fig:itsTable} shows the number of unpreconditioned MINRES iterations for each experiment with and without our inner-outer recycling approach. 

\subsection{Value of Inner-Outer Recycling} \label{subsec:numcomp}
There is a significant benefit to using two levels of recycling information.  
%To see this, we consider solving the full order model systems with unpreconditioned MINRES recycling as in \cite{KdS2006}. 
To see this, consider Algorithm \ref{alg:alg1} to construct the global basis.  We could solve the full order model systems in Step 1 (e.g. systems in \ref{eq:blkeqn}) with the unpreconditioned MINRES recycling approach in \cite{KdS2006}.  It is important to note that the recycle spaces would be different than those used in our new method.  Furthermore, in the new method we solve the correction equations (\ref{eq:resi}) as opposed to
solving (\ref{eq:blkeqn}).  For $j > 1$, the recycle spaces for the \cite{KdS2006} approach do not incorporate information from other systems corresponding to other right-hand sides.  In contrast, since we augment $\bfV$ from information about right-hand side $j$, we update $\bfK$.  The update in $\bfK$ then causes updates to $\bfr_{j+1}, \ldots, \bfr_{n_{src}+n_{det}}$, which are the right-hand sides in (\ref{eq:resi}).

Table~\ref{fig:tableRecyComp} compares the recycling of \cite{KdS2006} with our new approach.  The results show that with our approach, the number of iterations and the relative residuals decrease as you move from one right-hand-side to the next and also as you move from system to system.  The jump at right-hand-side $33$ is due to the fact that you are moving to the second half of the concatenated right-hand-sides, so these correspond to solving the adjoint problem.  Using the approach in \cite{KdS2006}, however, does not speed up convergence across right-hand-sides.   
%Updboth $\bfV$ and $\bfUU_j$ and subsection~\ref{subsec:cost} describes how this can be done efficiently.  
%In addition, using recycling you would need to still concatenate the solutions and do an SVD to get the reduced global basis in Algorithm \ref{alg:alg1}.  
In our approach, the reduced global basis, $\bfV$, is already constructed when we are done with the full order model solves.  We note there is a big difference in total number of MINRES iterations to squeeze all information from systems $1$ and $2$.  It took our approach $5,006$ iterations, while it took $22,659$ iterations for the recycling method in \cite{KdS2006}. 

\begin{table}[h]
        \centering
				\begin{tabular}{| c | c | c | c | c | c |}
				    \cline{1-6}
				      &  & \multicolumn{2} { c| } {Our Approach} & \multicolumn{2} { c| } {Recycling from \cite{KdS2006}} \\ \hline
				    System & RHS & Its & Initial Relative Residual & Its & Initial Relative Residual  \\ \hline
				    \multirow{6}{*}{1} & 1 & 140 & 7.523115e-05 & 140 & 7.523329e-05\\
				    
				    & 20 & 52 & 1.077645e-06 & 185 & 6.584893e-04 \\  
				    & 32 & 0 & 8.866398e-08 & 152 & 1.164802e-04 \\ 
				    & 33 & 124 & 4.776975e-05 & 127 & 5.213653e-05 \\ 
				    & 53 & 57 & 1.692149e-06 & 191 & 6.114295e-04 \\ 
				    & 64 & 0 & 9.960690e-08 & 151 & 1.450962e-04 \\ \hline
				    \multirow{6}{*}{2} & 1 & 118 & 4.673235e-05 & 131 & 5.091877e-05 \\ 
				    & 20 & 35 & 6.653493e-07 & 190 & 6.171251e-04 \\ 
				    & 32 & 0 & 8.754708e-08 & 153 & 1.051619e-04 \\ 
				    & 33 & 105 & 2.588292e-05 & 129 & 3.388392e-05 \\ 
				    & 53 & 37 & 8.454303e-07 & 188 & 6.405050e-04 \\  
				    & 64 & 0 & 7.960101e-08 & 151 & 9.958796e-05 \\ \hline
				    
				    \end{tabular}
        \caption{Comparison of MINRES recycling using $\bfUU_j$ as described in \cite{KdS2006} vs. the inner-outer approach using both  $\bfV$ and $\bfUU_j$ as described in Algorithm \ref{alg:recycbasis}.  Note that the two approaches lead to different choices for $\bfUU_j$ as well as different systems to solve.}
        \label{fig:tableRecyComp}
\end{table}

\section{Conclusions and Future Work} \label{sec:conclusions}

First, we established that our transfer function at zero frequency for DOT could be re-written in terms of a SPD matrix. Then, we developed an inner-outer Krylov recycling approach to update the global basis matrix relative to our new formulation of the transfer function.    
%By modifying traditional recycling with our inner-outer approach we are able to not only speed up the convergence of MINRES, but also build the global basis with only the non-redundant information.  
Two numerical experiments illustrate the success in using the ROM in place of the FOM during the optimization.  

In this paper, we only considered the $0$ frequency case.  It is non-trivial to extend the algorithm to the case when $\omega$ is non-zero, and it is therefore the subject of a forthcoming paper.  

Clearly, the performance of our method depends on the values of some parameters, such as the residual tolerance and the number of systems $K_*$.   We found, for example, that if we dropped the tolerance slightly, the number of system solves, and therefore the reduced model order, was even further reduced, without too much degradation in the reconstruction.  Likewise, using a larger value of $K_*$ gave slightly larger reduced order models, but with no improvement in the quality of the reconstruction.   The trade-offs in performance due to these selections are currently under investigation.   Finally, preliminary results indicate that solving the systems corresponding to different right-hand sides in a different ordering may also have an impact on the model order, and we will continue to 
investigate this phenomenon.  

%We are also investigating the order in which the full order model systems are solved.  Currently, we solve the full order model systems starting from the first right-hand-side and moving through until the last.  We are looking into reordering the right-hand-sides to improve convergence and our global basis.  

% -- THE END --

\vspace{6pt}
\noindent
%{\bf Acknowledgments.} We thank the anonymous reviewers {\em in advance} for their careful reading and helpful
%suggestions, which will greatly help us to improve this paper.

\bibliographystyle{abbrv}
\bibliography{npimr}

\end{document}